\documentclass[12pt]{article}
\usepackage[a4paper,headsep=0.6cm,left=2cm,right=2cm,top=2cm,bottom=2cm,twoside]{geometry}
\usepackage{graphicx,hyperref}

\usepackage{amsmath,amssymb,amsthm,bm}

\usepackage{amsfonts,euscript,enumitem}
\usepackage{graphicx}
\usepackage{mathrsfs}
\usepackage{tkz-graph}
\usetikzlibrary{decorations.markings}
\usepackage{color,pgf,tikz,pgfbaseimage,tikz-cd}
\usetikzlibrary{arrows,cd,calc}

\title{On $\delta$-deformations of polygonal dendrites.}

\author{Dmitry Drozdov \and Mary Samuel \and Andrey Tetenov \footnote{Supported by Russian Foundation of Basic Research projects 16-01-00414  and 18-501-51021}}
\begin{document}
\sloppy

\renewcommand{\refname}{References}
\renewcommand{\proofname}{Proof.}
\renewcommand{\figurename}{Fig.}

\newcommand{\rr}{\mathbb{R}}
\newcommand {\nn} {\mathbb{N}}
\newcommand {\zz} {\mathbb{Z}}
\newcommand {\bbc} {\mathbb{C}}
\newcommand {\rd} {\mathbb{R}^d}
\newcommand {\rpo}{\mathbb{R}_+^1}

 \newcommand {\al} {\alpha}
\newcommand {\be} {\beta}
\newcommand {\da} {\delta}
\newcommand {\Da} {\Delta}
\newcommand {\ga} {\gamma}
\newcommand {\Ga} {\Gamma}
\newcommand {\la} {\lambda}
\newcommand {\La} {\Lambda}
\newcommand{\om}{\omega}
\newcommand{\Om}{\Omega}
\newcommand {\sa} {\sigma}
\newcommand {\Sa} {\Sigma}
\newcommand {\te} {\theta}
\newcommand {\vte} {\vartheta}
\newcommand {\fy} {\varphi}
\newcommand {\Fy} {\varPhi}
\newcommand{\ep}{\varepsilon}
\newcommand{\ro}{\varrho}

\newcommand {\ra} {\rightarrow}
\newcommand{\IN}{{\subset}}
\newcommand{\NI}{{\supset}}
\newcommand {\mmm}{{\setminus}}
\newcommand{\8}{{\infty}}
\newcommand{\io}{{I^\infty}}
\newcommand{\ia}{{I^*}}
\newcommand{\0}{{\varnothing}}
\newcommand{\vse}{$\blacksquare$}
\def \Ln {\mathop{\rm Ln}\nolimits}

\newcommand {\Dl} {\Delta}

\newcommand{\bj}{{\bf {j}}}
\newcommand{\bi}{{\bf {i}}}
\newcommand{\bk}{{\bf {k}}}
\newcommand{\bl}{{\bf {l}}}
\newcommand{\bu}{{\bf {u}}}

\newcommand{\eJ}{{\EuScript J}}
\newcommand{\eC}{{\EuScript C}}
\newcommand{\wP}{{\widetilde P}}
\newcommand{\eU}{{\EuScript U}}
\newcommand{\eP}{{\EuScript P}}
\newcommand{\eS}{{\EuScript S}}
\newcommand{\eW}{{\EuScript W}}
\newcommand{\eV}{{\mathcal V}}
\newcommand{\eZ}{{\EuScript Z}}
\newcommand{\eK}{{\EuScript K}}

\newcommand{\beq}{\begin{equation}}
\newcommand{\eeq}{\end{equation}}

\def \Id {\mathop{\rm Id}\nolimits}
\def \diam {\mathop{\rm diam}\nolimits}
\def \sup {\mathop{\rm sup}\nolimits}
\def \fix {\mathop{\rm fix}\nolimits}
\def \Lip {\mathop{\rm Lip}\nolimits}
\def \min {\mathop{\rm min}\nolimits}

\newtheorem{thm}{\bf Theorem}
 \newtheorem{cor}[thm]{\bf Corollary}
 \newtheorem{lem}[thm]{\bf Lemma}
 \newtheorem{prop}[thm]{\bf Proposition}
 \newtheorem{dfn}[thm]{\bf Definition}
 \theoremstyle{definition}
\newtheorem{rmk}{Remark}

\newcommand{\dok}{{\bf{Proof:  }}}
\newcommand{\red}{\textcolor{red}}
\newcommand{\yellow}{\textcolor{yellow}}
\newcommand{\blue}{\textcolor{blue}}
\newcommand{\green}{\textcolor{green!60!black}}

\maketitle

\begin{abstract}
We find the conditions under which the attractor $K(\eS')$ of a deformation $\eS'$ of a contractible polygonal system $\eS$ is a dendrite.
\end{abstract}


\bigskip

This is a very convenient though rather restrictive way to define post-critically finite self-similar dendrites in the plane using contractible P-polygonal systems. This approach was discussed in \cite{TSV0},\cite{TSV},\cite{TSM}. It turned out that well-known examples of self-similar dendrites are obtained using such systems. Nevertheless, if we move slightly the   vertices of the main polygon $P$ and of polygons $P_i$, defining the polygonal system $\eS$, and change the system $\eS$ accordingly, we often obtain a system $\eS'$ of a more general type whose attractor $K'$ is a dendrite too. We call such systems generalized polygonal systems (Definition \ref{gps}) and in the case when polygons $P_i'$ differ from the polygons $P_i$ less than by $\da$, we call such systems  $\da$-deformations (Definition \ref{deform}) of the polygonal system $\eS$. In this paper we begin initial study of  generalized polygonal systems and  $\da$-deformations of contractible polygonal systems.

In Theorem \ref{pcint} we formulate sufficient conditions under which the attractor $K$ of a generalized polygonal system $\eS$ is a dendrite. These conditions are expressed in terms of intersections $K_i\cap K_j$ of the pieces of the attractor $K$. 
In Theorem \ref{attrmap} we show that a $\da$-deformation $\eS'$ of a contractible polygonal system $\eS$  defines a continuous map $\hat f: K\to K'$ of  respective attractors of these systems which agrees with the action of $\eS$ and $\eS'$ and give conditions under which $\hat f$ is a homeomorphism.
In Theorem \ref{PMT} we show that Parameter Matching Condition is a necessary condition for a generalized polygonal system to generate a dendrite.
In Theorem \ref{mainthm} we show that if $\da$ is sufficiently small and the system $\eS'$ is  $\da$-deformation of a contractible P-polygonal system $\eS$, which satisfies Parameter Matching Condition, then the attractor $K(\eS')$ is a dendrite, homeomorphic to $K(\eS)$.
 
\subsection{Preliminaries}
\subsubsection{Self-similar sets}
\begin{dfn} 
Let $\eS=\{S_1, S_2, \ldots, S_m\}$ be a system of (injective) contraction maps on the complete metric space $(X, d)$.
 A nonempty compact set $K\IN X$ is called the attractor of the system $\eS$, if $K = \bigcup \limits_{i = 1}^m S_i (K)$. \end{dfn}
 
  The system $\eS$ defines its Hutchinson operator $T$  by  $T(A) = \bigcup \limits_{i = 1}^m S_i (A)$. By Hutchinson's Theorem, the attractor $K$ is unique for $\eS$ and for any compact set $A\IN X$ the sequence $T^n(A)$ converges to $K$. We also call the   subset $K \IN X$ self-similar with respect to $\eS$.\\
Throughout the whole paper, the maps $S_i\in \eS$ are supposed to be  similarities and the set $X$ to be $\mathbb{R}^2$. We will use complex notation for the point on the plane so each 
similarity  will be written as $S_j(z)=q_je^{i\al_j}(z-z_j)+z_j$, where $q_j=\Lip S_j$ and $z_j=\fix(S_j)$. For a system $\eS$, let $q_{min}=\min\{q_j,j\in I\}$ and $q_{max}=\max\{q_j,j\in I\}$.\\
Here $I=\{1,2,...,m\}$ is the set of indices, while $\ia=\bigcup\limits_{n=1}^\8 I^n$ 
is the set of all finite $I$-tuples, or multiindices $\bj=j_1j_2...j_n$. The length $n$ of the multiindex $\bj=j_1...j_n$ is denoted by  $|\bj|$ and $\bi\bj$  denote the concatenation of the corresponding multiindices. We say $\bi\sqsubset\bj$, if  $\bj=\bi\bl$ for some $\bl\in\ia$; if $\bi\not\sqsubset\bj$ and $\bj\not\sqsubset\bi$, $\bi$ and $\bj$ are {\em incomparable}.\\ 
 
For a multiindex $\bj\in\ia$ we write
$S_\bj=S_{j_1j_2...j_n}=S_{j_1}S_{j_2}...S_{j_n}$ 
and for the set $A
\subset X$ we denote $S_\bj(A)$ by $A_\bj$; 
we also denote by $G_\eS=\{S_\bj, \bj\in\ia\}$ the semigroup, generated by $\eS$;\\
$I^{\8}=\{{\bf \al}=\al_1\al_2\ldots,\ \ \al_i\in I\}$ denotes the
index space; and $\pi:I^{\8}\rightarrow K$ is the {\em index map
 }, which sends $\bf\al$ to  the point $\bigcap\limits_{n=1}^\8 K_{\al_1\ldots\al_n}$.\\
 Along with a system $\eS$ we will consider its n-th refinement $\eS^{(n)}=\{S_\bj, \bj\in I^n\}$, whose Hutchinson's operator is equal to $T^n$.

\begin{dfn}The system ${\eS}$ satisfies the {\em open set condition} (OSC) if there exists a non-empty open set $O \IN X$ such that $S_i (O), \{1 \le i\le m\}$ are pairwise disjoint and all contained in $O$.\end{dfn}

Let $\eC$ be the union of all  $S_i(K)\cap S_j(K)$, $i,j \in I, i\neq j$. {\em The post-critical set} $\eP$ of the system $\eS$ is the set of all $\alpha\in I^{\8}$ such that for some ${\bf j}\in I^*$, $S_ {\bf j}(\alpha)\in\eC$. In other words, $\eP= \lbrace \sigma^k(\alpha) : \alpha\in \eC,k\in \mathbb{N}\rbrace$, where the map $\sigma^k:I^{\8}\to I^{\8}$ is defined by 
$\sigma^k(\al_1\al_2\ldots)=\al_{k+1}\al_{k+2}\ldots$
A system $\eS$ is called {\em post-critically finite} (PCF) \cite{Kig} if its post-critical set $\eP$ is finite. Thus, if the system  $\eS$ is postcritically finite then there is a finite set $\eV=\pi(\eP)$ such that for any non-comparable $\bi,\bj\in\ia$,  $K_\bi\cap K_\bj= S_\bi(\eV)\cap S_\bj(\eV)$.

\subsubsection{Dendrites}
A {\em dendrite} is a locally connected continuum con\-taining no simple closed curve. 

The order $Ord(p,X)$ of the point  $p$  with respect to a dendrite  $X$ is  the number  of components of the set $X \setminus \{p\}$.  
Points of order 1 in a dendrite $X$ are called {\em end points} of $X$;   a point $p\in X$ is called a {\em cut point} of $X$ if $X \setminus \{p\}$ is disconnected; 
points of order at least 3 are called {\em ramification points} of $X$.

A continuum $X$  is a dendrite iff
X is locally connected and uniquely arcwise connected.

\subsubsection{Contractible polygonal systems}
 Let $P \subset {\mathbb{R}}^2$ be a finite polygon homeomorphic to a disk, $ \eV_P= \{A_1, \ldots, A_{n_P}\}$ be the set of its vertices. Let also $\Om(P,A)$ denote the angle with vertex $A$ in the polygon $P$. We 
consider  a system of similarities $\mathcal{S} = \{S_1, \ldots, S_m\}$ in ${\mathbb{R}}^2$ such that:\\
{\bf (D1)} for any $i \in I$ set $P_i = S_i (P) \subset P$;\\ 
{\bf (D2)} for any $i \ne j, i, j \in I, P_i \bigcap P_j =  \eV_{P_i} \bigcap  \eV_{P_j}$ and $\#(\eV_{P_i} \bigcap  \eV_{P_j})<2$;\\  
$\bf (D3)$ $\eV_P \subset \bigcup \limits_{i \in I} S_i (\eV_P)$;\\
\\
 $\bf (D4)$ the set $\widetilde P = \bigcup \limits_{i = 1} ^m P_i$ is contractible.\\
 \begin{dfn}\label{cps} The system  $\eS$ satisfying the conditions $
\bf (D1 - D4)$, is called a {\it contractible $P$-polygonal system of similarities}.
\end{dfn}

This theorem was proved by the authors in \cite[Theorem 4,(g)]{TSV}(or \cite[Theorem 10,(g)]{TSVa}):
\begin{thm}\label{refsys} Let  $\eS$ be a contractible $P$-polygonal  system  of similarities.\\
(a) The system $\eS$ satisfies (OSC).\\
(b) $P_\bj\IN P_\bi$ iff  $\bj\sqsupset\bi$.\\
(c) If $\bi\sqsubset\bj$, then  $ S_\bi(\eV_P)\cap P_\bj\IN S_\bj(\eV_P)$.\\
(d) For any incomparable $\bi,\bj\in \ia$, $\#(P_\bi\cap P_\bj)\le 1$ and $P_\bi\cap P_\bj=S_\bi(\eV_P)\cap S_\bj(\eV_P)$.\\
(e) The set $G_\eS(\eV_P)$ of vertices of polyhedra $P_\bj$ is contained in $K$.\\
(f)If $x\in K\mmm G_\eS(\eV_P)$, then $\#\pi^{-1}(x)=1$.\\
(g) For any $x\in G_\eS(\eV_P)$ there is   $\ep>0$ 
and a finite system $\{\Om_1,...,\Om_n\},$ where $n=\#\pi^{-1}(x)$,
of disjoint  angles with vertex  $x$,
such that if $x\in P_\bj$ and $\diam P_\bj<\ep$, then for some
  $k\le n$, $\Om(P_\bj,x)=\Om_k$. Conversely, for any $\Om_k$ there is such $\bj\in\ia$, that  $\Om(P_\bj,x)=\Om_k$.
\end{thm}
\resizebox{.9    \textwidth}{!}{{
\begin{tikzpicture}[line cap=round,line join=round,>=stealth ,scale=1]\scriptsize
\node at (-6.3,-2.3) {
   \includegraphics[width=.57\textwidth]{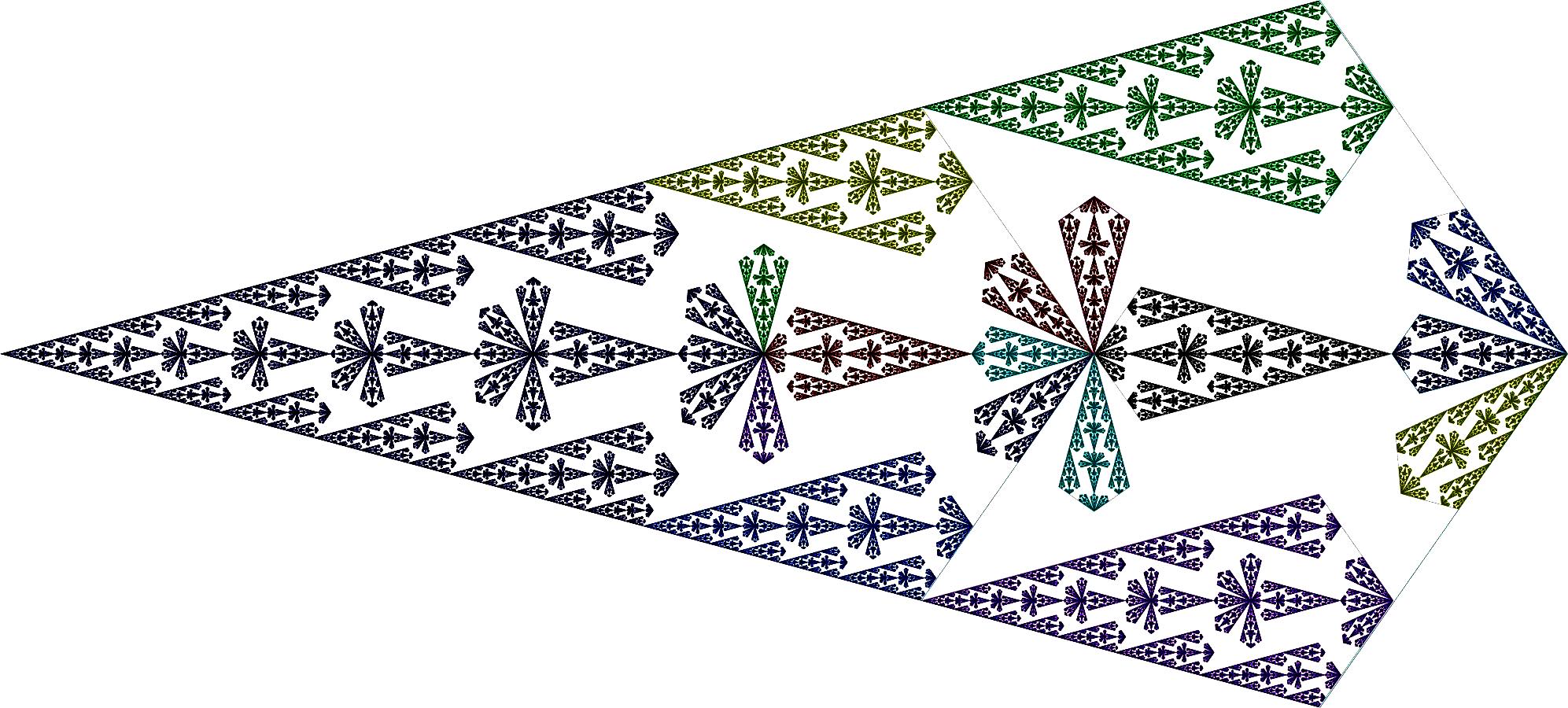}}; \node at (2,-2) {\includegraphics[width=.35\textwidth]{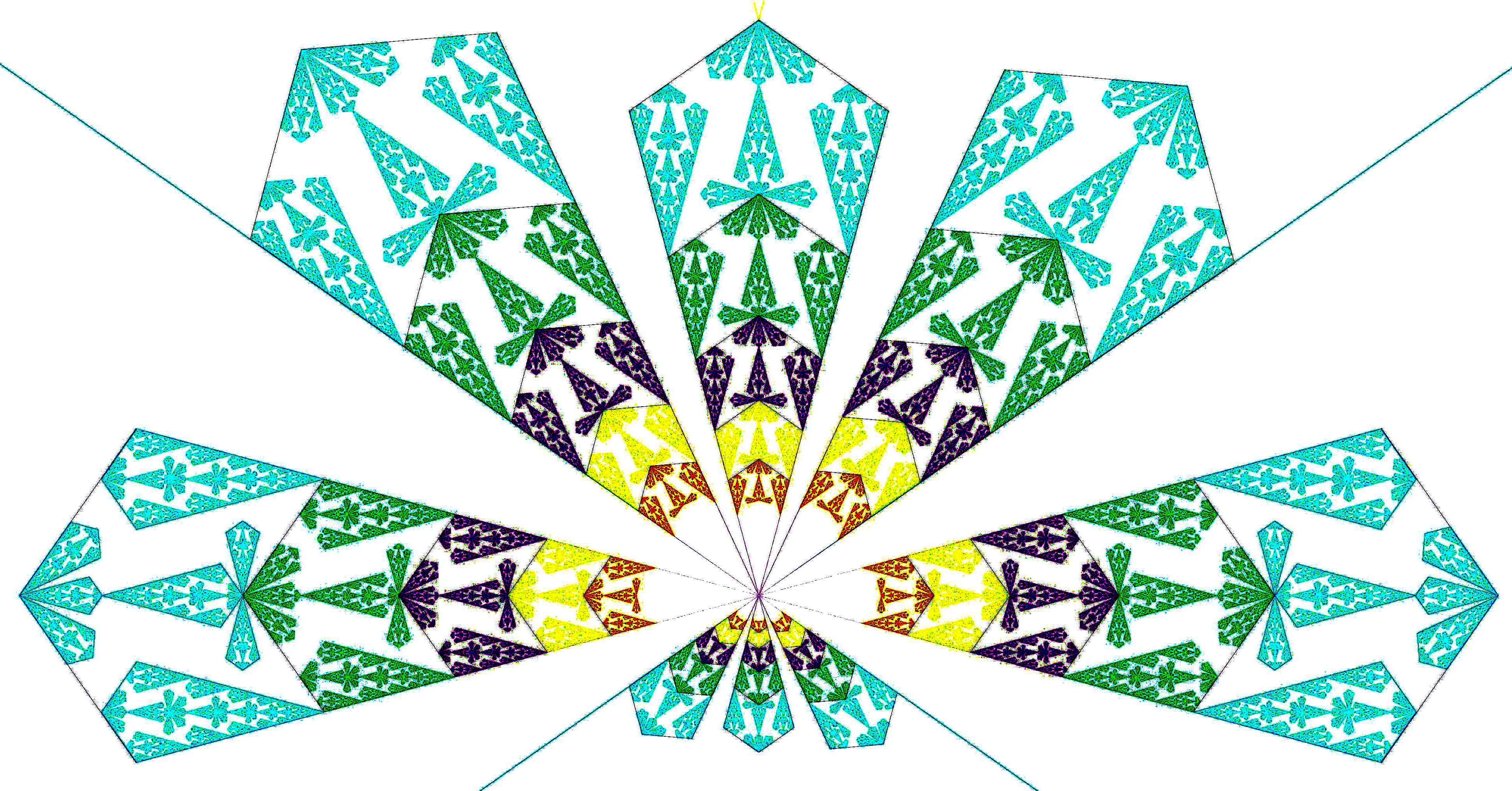}};

   \node at(-5.3,-.25){$B$}; \node at (-8.5,-4) {\small Polygonal system and its attractor};\node at (2,-4) {\small Local structure of  $K$ near the vertex B.(rotated)};\node at (-4.15,-3.35) {$\scriptsize{\rho_0}$};
\path[thick,red]
	   (-5.3,-.4) edge[->]  (-4.44,-2.2);
	   \path[thin,red] (-4.38,-3.27) edge[<->]  (-4.3,-3.5);
\end{tikzpicture}}} 

\begin{dfn}\label{cyclic}
Let $\eS$ be a contractible $P$-polygonal system of similarities. The vertex $A \subset \eV_P$ is called a cyclic vertex, if there is such multiindex  $\bi=i_1 i_2 \ldots i_k$, that $S_\bi(A)=A$. The least number $k=|\bi|$ among all $\bi$ for which $S_\bi(A)=A$ is called {\em the order} of the cyclic vertex $A$.
\end{dfn}

\begin{dfn}
A point $B \in\cup_{i=1}^m  \eV_{P_i}$ is subordinate to a cyclic vertex  $A$, if for certain multiindex $\bi, S_{\bi}(A)=B$.
\end{dfn}

\begin{prop}\label{1ordcyc}
Let $\eS$ be a contractible $P$-polygonal system of similarities. Then:\\
(1) Each vertex $B\in \eV_P$ is subordinate to some cyclic vertex.\\
(2) There is such $n$, that in the system $\eS^{(n)}=\{S_\bi, \bi\in I^n\}$ all the cyclic vertices have order 1.
\end{prop}
{\bf Proof.} Notice that if $A\in \eV_P$ is a cyclic vertex, then there is such $\bj\in\ia$ that $S_\bj(A)=A$. Therefore if for some $\bj\in\ia$, $ A\in P_\bj$, then for some $n$, $S_\bj^n(P)\IN P_\bj \IN P$, $A$ being a vertex of each of these polygons. Since $\Om(S_\bj^n(P),A)=\Om(P,A)$, for any $\bj\in\ia$, for which $A\in P_\bj$, $\Om(P_\bj,A)=\Om(P,A)$. This implies that $\#\pi^{-1}(A)=1$ and for any $n$ there is unique $\bj\in I^n$ such that $A\in P_\bj$. 

Conversely if for any $\bi\in\ia$, for which $A\in P_\bi$, $\Om(P_\bi,A)=\Om(P,A)$ then $\#\pi^{-1}(A)=1$ and $A$ is a cyclic vertex of the system $\eS$.

Then, by Theorem \ref{refsys}, for any vertex $B\in G_\eS(\eV_P)$ there is a finite set $\{\bi_1,...,\bi_n\}$ of incomparable multiindices  such that for any $l,l'$, $P_{\bi_l}\cap P_{\bi_{l'}}=\{B\}$, the set $\bigcup\limits_{l=1}^k K_{\bi_l}$ is a neighbourhood of the point $B$ in $K$ and for any $l=1,...,k$, the point $S_{\bi_l}^{-1}(B)=A_l$ is a cyclic vertex. This proves (1).

Let now $A_1,....,A_k$ be the full set of cyclic vertices in $\eV_P$ and $p_1,...,p_k$ be their respective orders. Let $N$ be the  least common multiple of $p_1,...,p_k$. Then $\eS^{(n)}$ is the desired $P$-polygonal system.\vse

\subsubsection{Main parameters of a contractible polygonal system}

For any set $X\IN\rr^2$ or point $A$ by $V_\ep(X)$ (resp.$V_\ep(A)$) we  denote $\ep$-neighbourhood of the set $X$ (resp. of the point $A$) in the plane.\\

$\bm{\rho_0:}$  Take such $\rho_0>0$ that:\\ (i) for any vertex $A\in \eV_P$,   $V_\rho(A) \bigcap P_k \ne \0 \Rightarrow A\in P_k$;\\ (ii) for any $x, y \in P$ such that there are $P_k, P_l: x \in P_k, y \in P_l$ and $P_k \bigcap P_l = \0, d (x, y) \ge \rho_0$. \\ 

\resizebox{.9    \textwidth}{!}{{
\begin{tikzpicture}[line cap=round,line join=round,>=stealth ,scale=1]\scriptsize
\node at (-6.3,-2) {
   \includegraphics[width=.57\textwidth]{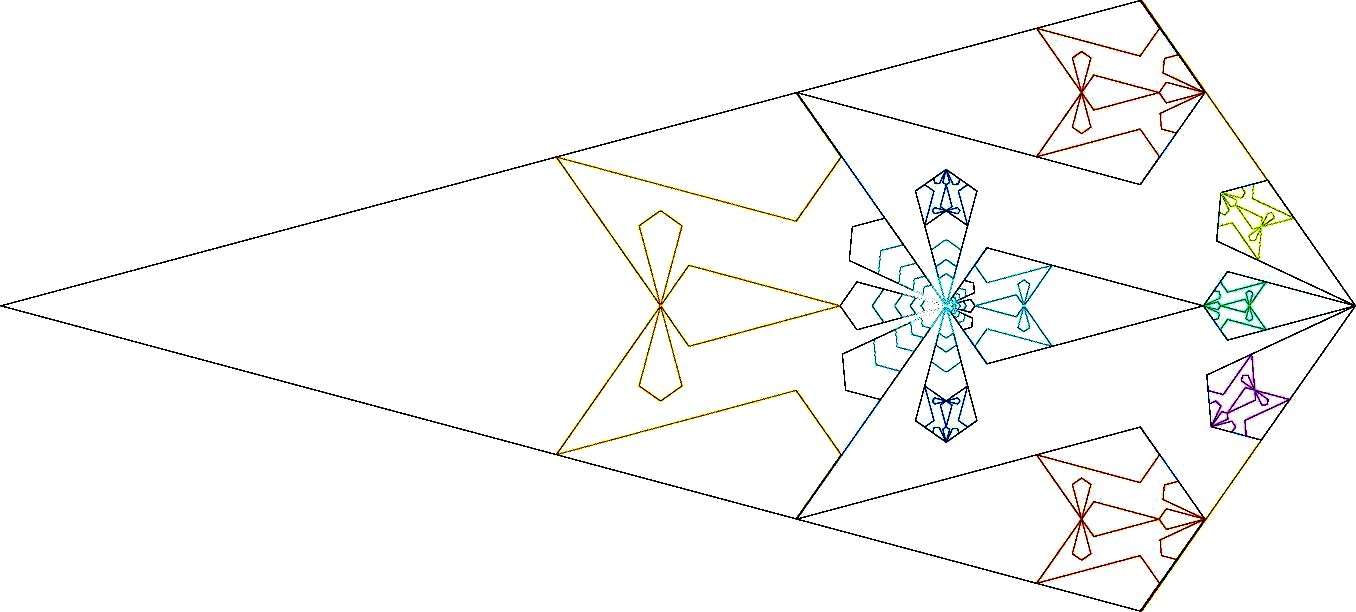}}; \node at (2,-2.68) {\includegraphics[width=.35\textwidth]{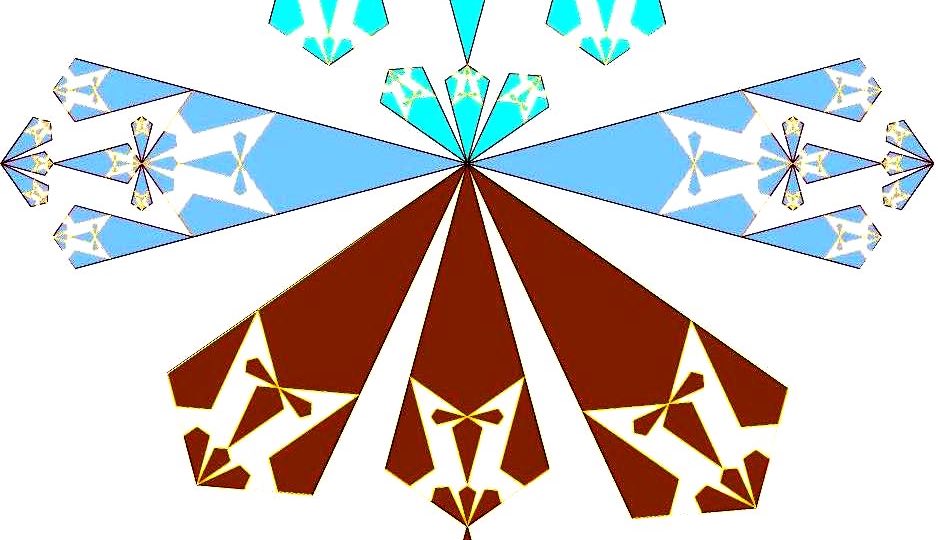}};

  	  \path[thick,red] (2,-2) edge[->]  (4.7,-3.2);
	  \path[thin,red] (-1.46,-2) edge[]  (-2.8,-2.62);
	  \path[thin,red] (-1.46,-2) edge[]  (-2.9,-2.39);
	    \draw[thick,red] (4.1,-4.1) arc (-45:20:2.96);
	     \draw[thick,red] (-.1,-4.1) arc (-135:-200:2.96);
	     \draw[thin,red] (-2.55,-2.5) arc  (-155:-165.5:1.08);
	    \draw [thick,blue](2,-2) circle [radius=.35];
	    \draw [thick,red](-4.38,-2) circle [radius=0.98];
	    \draw [thin,blue](-4.38,-2) circle [radius=0.12];
	    \node at(4.4,-3.3){$\bm\rho_2$};\node at(1.5,-2.2){$\bm\rho_1$};\node at(-2.8,-2.5){$\bm\al_0$};
	    \node at (-3,-4.5) {\small Choosing the parameters $\al_0$, $\rho_1$ and $\rho_2$ for a polygonal system.};
\end{tikzpicture}}}

 $\bm{\rho_1,\rho_2:}$    As it follows from Theorem \ref{refsys}, for any vertex $B\in \eV_\wP$  there is a finite set of cyclic vertices $A_{i_1},...,A_{i_k}\in \eV_P$,  and multiindices
$\bj_1,...,\bj_k$  such that for any $l=1,...,k$, $S_{\bj_l}(A_l)=B$  and $S_{i_l}(A_l)=A_l$  and the set $\bigcup\limits_{l=1}^k S_{\bj_l}S_{i_l}^n(K)$ is a neighborhood of  the point $B$ in $K$ for any $n\ge 0$.\\
Let $\rho_1$ and $\rho_2$ be such positive numbers that for for any vertex $B\in \eV_\wP$  \beq\label{ro12} (V_{\rho_1}(B)\cap K)\IN \bigcup\limits_{l=1}^k S_{\bj_l}(P_{i_l})\mbox{\qquad   and  \qquad  }
  \bigcup\limits_{l=1}^k P_{\bj_l}\IN V_{\rho_2}(B).\eeq
  
$\bm{\al_0:}$  Let $\al_0$ denote the minimal angle between those sides of polygons $P_i, P_j, i,j\in I$, which have common vertex.

  {\bf Arrangement of maps fixing cyclic vertices.} Let $\eS$ be a contractible $P$-polygonal system all of whose cyclic vertices have order 1. In this case we can arrange the indices in $I$ and enumerate the vertices in $\eV_P$ in such way that  each cyclic vertex $A_l$ will be the fixed point of $S_l\in \eS$. Notice that $S_l$ is a homothety $S_l(z)=q_l(z-A_l)+A_l$ and the polygon $P$  lies inside the angle $\Om(P, A_l)$ and $K\mmm \{A_l\} =\bigsqcup\limits_{n=0}^\8S_l^n(K\mmm K_l)$. The number of points in $\overline{K_l\mmm S_l(K_l)}\cap S_l(K_l)$  is finite and is equal to the ramification order of $A_l$ in $K$.

\subsection{Generalized polygonal systems.} 

If we omit the condition {\bf(D1)} in the definition of contractible $P$-polygonal system $\eS$, we get the definition of a {\em generalized $P$-polygonal system}:
\begin{dfn}\label{gps}
A system \ $\eS=\{S_1,...,S_m\} $, satisfying the conditions {\bf D2-D4},
is called a generalized $P$-polygonal system of similarities.
\end{dfn}

 \hspace{-1em}\resizebox{.9    \textwidth}{!}{{\small\begin{tikzpicture}[line cap=round,line join=round,>=stealth ,x=10.0cm,y=10.0cm]
    \node  at (-.098,0.02) {\includegraphics[width=.36\textwidth]{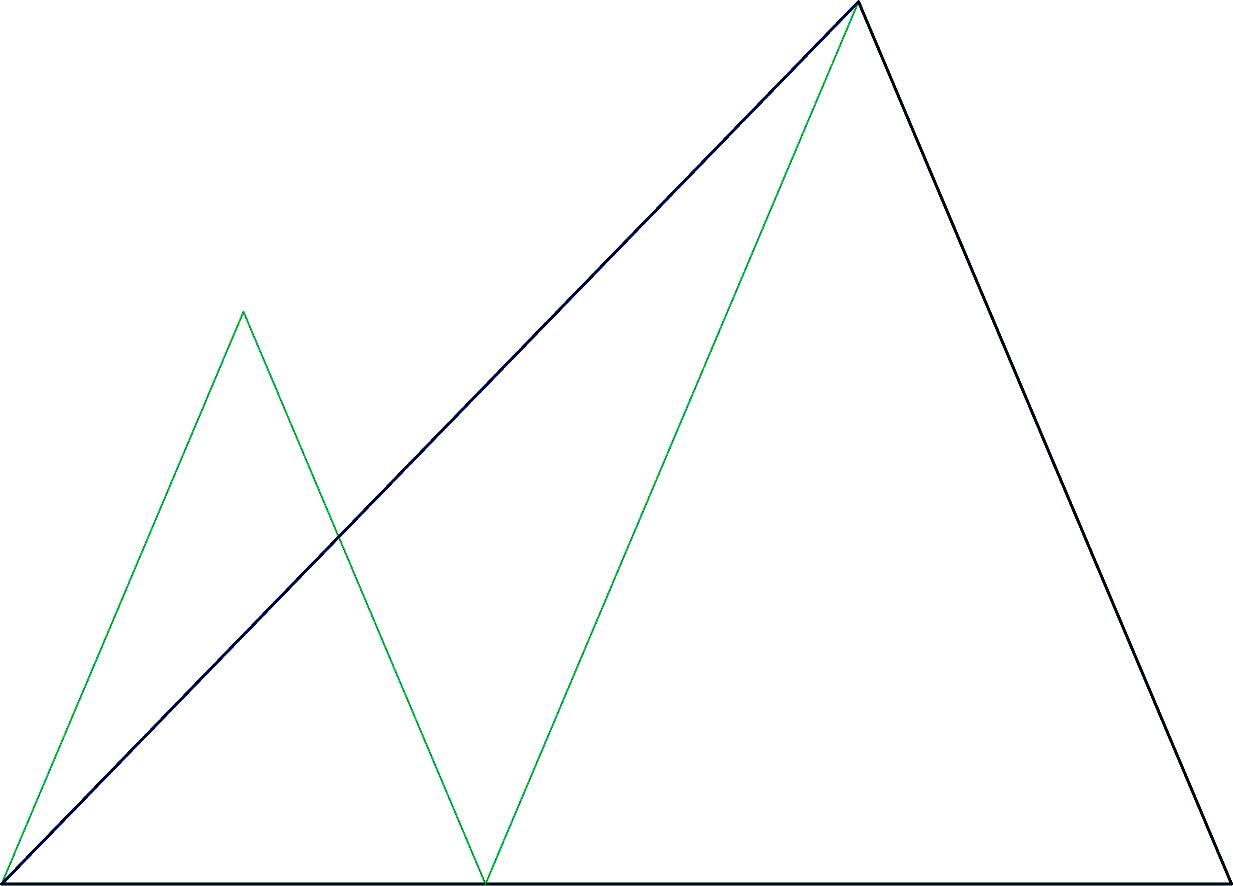} };
   \node  at (.65,-.01) {\includegraphics[width=.37\textwidth]{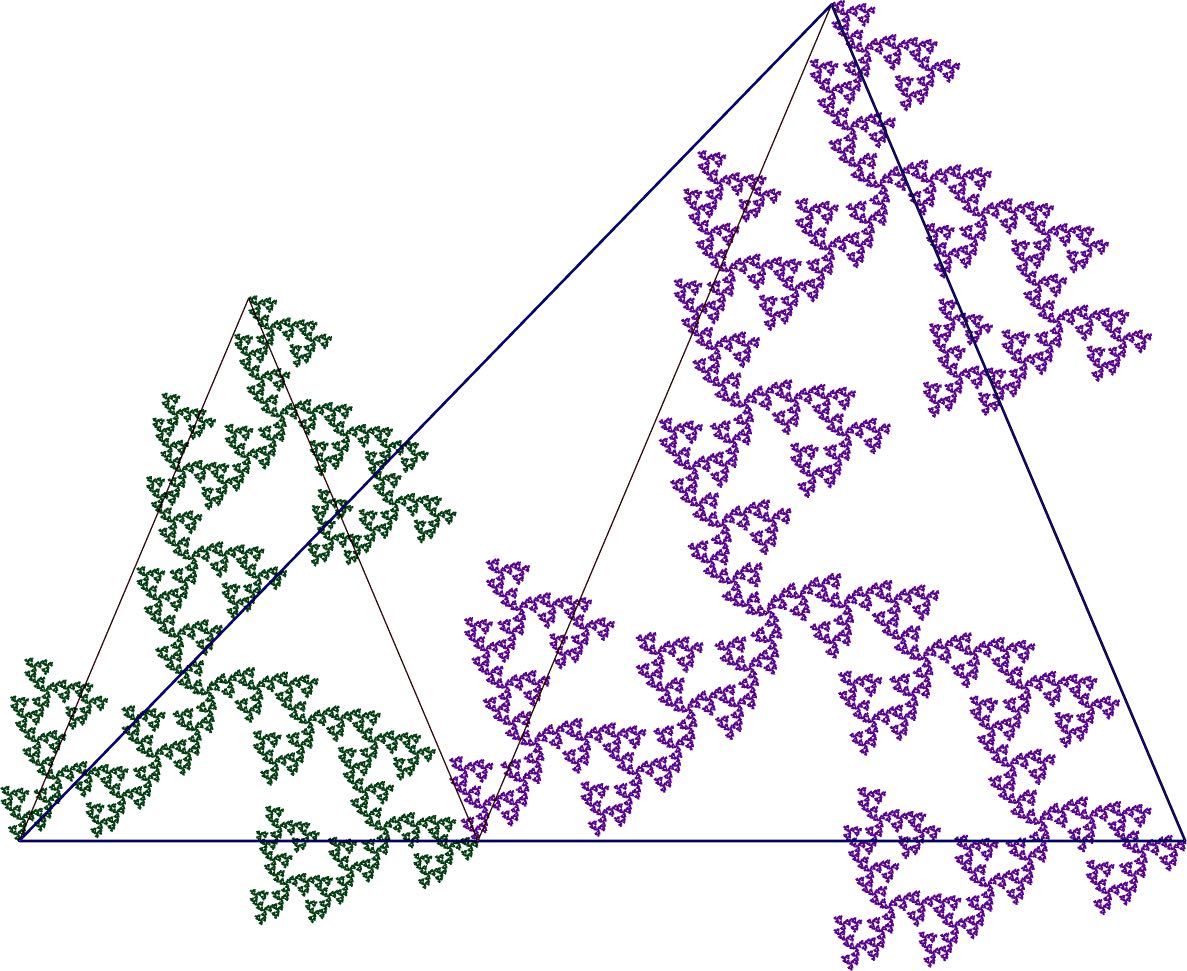}};
      \node at (0.02,0){$P_2$}; \node at (-0.14,-.02){$P$};\node at (-0.28,-0) {$P_1$};\node at (0.51,-.08) {$K_1$};\node at (.81,-.03){$K_2$};
\end{tikzpicture}}}

\begin{thm}\label{pcint}
Let $\eS$ be a generalized $P$-polygonal system.  If for any $i, j \in I$ \begin{equation}\label{icnd}S_i(K)\cap S_j(K)=P_i\cap P_j,\end{equation}
then the attractor $K$ of the system $\eS$ is a dendrite.
\end{thm}

\dok  Let $i,i'\in I$. By a (simple) chain of indices, connecting $i$ and $i'$, we mean a sequence $i=i_1,i_2,...,i_l=i'$ of pairwise different indices such that $P_{i_k}\cap P_{i_{k'}}=\0$  if $|k'-k|>1$, and that for any $k=1,...,l-1,$\qquad $P_{i_k}\cap P_{i_{k+1}}=\{x_k\}$, where $x_k$ denotes a common vertex of the polygons $P_{i_k}$ and $P_{i_{k+1}}$. The last condition also means, that  $K_{i_k}\cap K_{i_{k+1}}\ni x_k$ for any $k\in I$.\\ 

Since in a generalized polygonal system for any two indices $i,i'$ there is a chain of indices $i=i_1,i_2,...,i_l=i'$ connecting them, then by \cite[Theorem 1.6.2]{Kig}, the attractor $K$ is connected, locally connected  and arcwise connected. Thus, any two points of $K$ can be connected by some Jordan arc in $K$.\\

 Notice also that  if the condition
 (\ref{icnd}) holds, and the indices $i,i'\in I$ can be connected by a chain $i=i_1,i_2,...,i_l=i'$, then for any points $x\in K_i$, $y\in K_{i'}$ there is some Jordan arc $\ga_{xy}\in K$, consisting of subarcs \begin{equation}\label{subarcs}\ga_{xx_1}\IN K_{i_1},\ldots,\ \ga_{x_{k-1}x_k}\IN K_{i_k},\ldots,\ \ga_{x_{l-1}y}\IN K_{i_l}\end{equation} with disjoint interiors.\\
 
   At the same time, if the condition
 (\ref{icnd}) holds, and a Jordan arc  $\ga_{xy}\IN K$ with endpoints in $x$ and  $y$,  meets sequentially the pieces $K_{i_1},...,K_{i_l}$,  then it passes sequentially through  the points   $x_k$, where $\{x_k\}=K_{i_{k-1}}\cap K_{i_{k}}$ and consists of  subarcs of the form (\ref{subarcs}) with disjoint interiors.\\ 
 
 And vice versa,  if the condition (\ref{icnd}) holds, then for any Jordan arc $\ga_{xy}$  in $K$ there is unique chain of indices  $i_1,\ldots,i_l$, such that $\ga_{xy}$ consists of  subarcs of the form (\ref{subarcs}).\\
 
  We need a small Lemma to continue the proof:

 \begin{lem}\label{arcs} Let $\bj\in I^*$ and $x,y\in K_\bj$.
 If the condition (\ref{icnd}) holds, then for any two Jordan arcs   $\la_1,\la_2$ with endpoints $x,y$, the distance $d_H(\la_1,\la_2)\le q_{max}\diam K_\bj$. \end{lem}
 
 \dok Indeed, consider the Jordan arcs $\la_1'=S_\bj^{-1}(\la_1)$ and $\la_2'=S_\bj^{-1}(\la_2)$, connecting $x'=S_\bj^{-1}(x)$ and $y'=S_\bj^{-1}(y)$ in $K$. Let
 $x'\in K_i$ and $y'\in K_{i'}$, and let $i_1,i_2,...,i_l$ be the chain, connecting $i$ and $i'$. Then each of the arcs $\la_1',\la_2'$ consists of subarcs, connecting sequentially the pairs of points $x_k,x_{k+1}$ in the sequence $x', x_1,...,x_{l-1},y'$, and lying in respective pieces $K_{i_k}$.
 Since the diameters of these sets are not greater than $q_{max}\diam K$, $d_H(\la'_1,\la'_2)\le q_{\max}\diam K$. Then $d_H(\la_1,\la_2)\le q_{\max}\diam K_\bj\le \diam K q_{\max}^{|\bj|+1}$.\vse\\
 
 Now we can finish the proof of the Theorem. Let $\la$ and $\la'$ be Jordan arcs in $K$ with endpoints at $x$ and $y$. Applying the Lemma  \ref{arcs} by induction to the subarcs of which the arcs $\la$ and $\la'$ consist, we get that for any $n>|\bj|$, $d_H(\la_1,\la_2)\le q_{\max}^n\diam K$. Taking a limit with $n\to\8$, we obtain that a Jordan arc, connecting the points $x$ and $y$ is unique. Therefore  $K$ is a dendrite.\vse\\

\begin{rmk} It is possible for a generalized $P$-polygonal system $\eS$ not to satisfy the condition \ref{icnd} and to have the attractor $K$ which is a dendrite. The attractor $K$ of a generalized polygonal system $\eS$ on the picture below
  is a dendrite, but  $P_7\cap P_9=\0$, whereas $K_7\cap K_9$ is  a line segment.\end{rmk}
  
  \hspace{-1em}\resizebox{.9    \textwidth}{!}{{\begin{tikzpicture}[line cap=round,line join=round,>=stealth ,x=10.0cm,y=10.0cm]
\node[anchor=south west,inner sep=0] at (0,-.227) {
   \includegraphics[width=.45\textwidth]{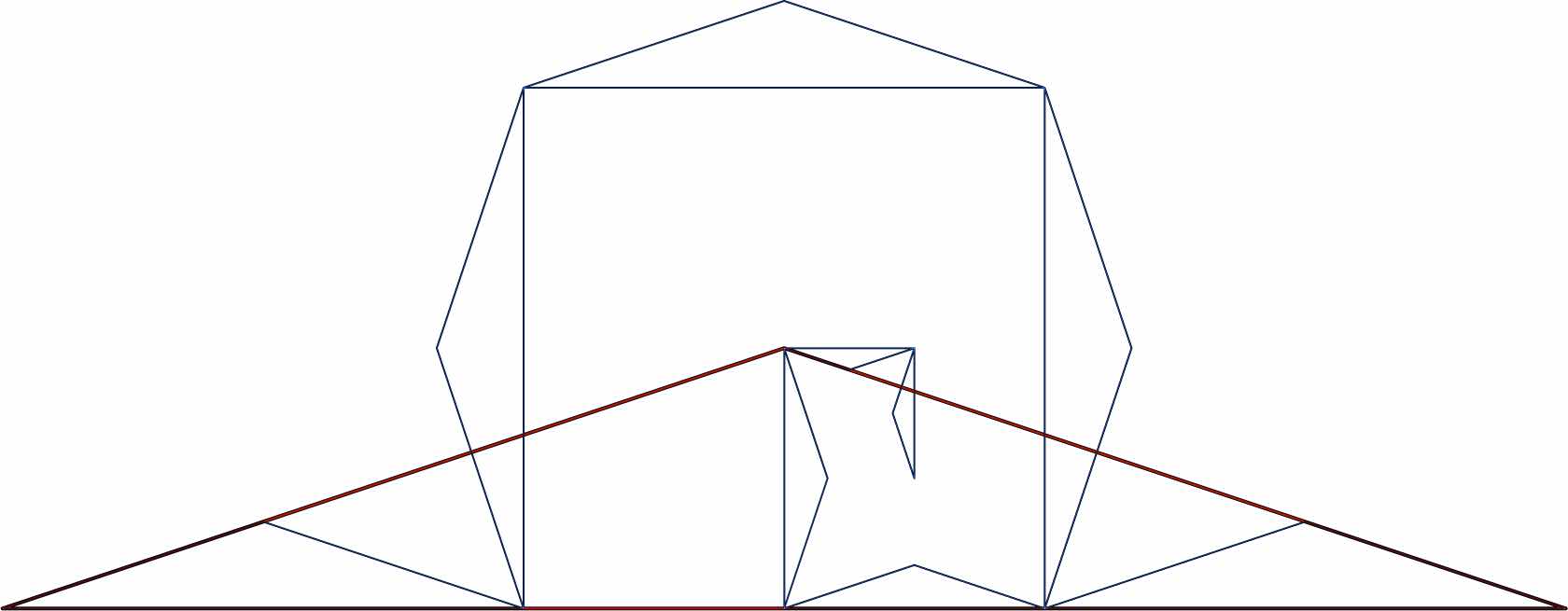}\quad  \includegraphics[width=.45\textwidth]{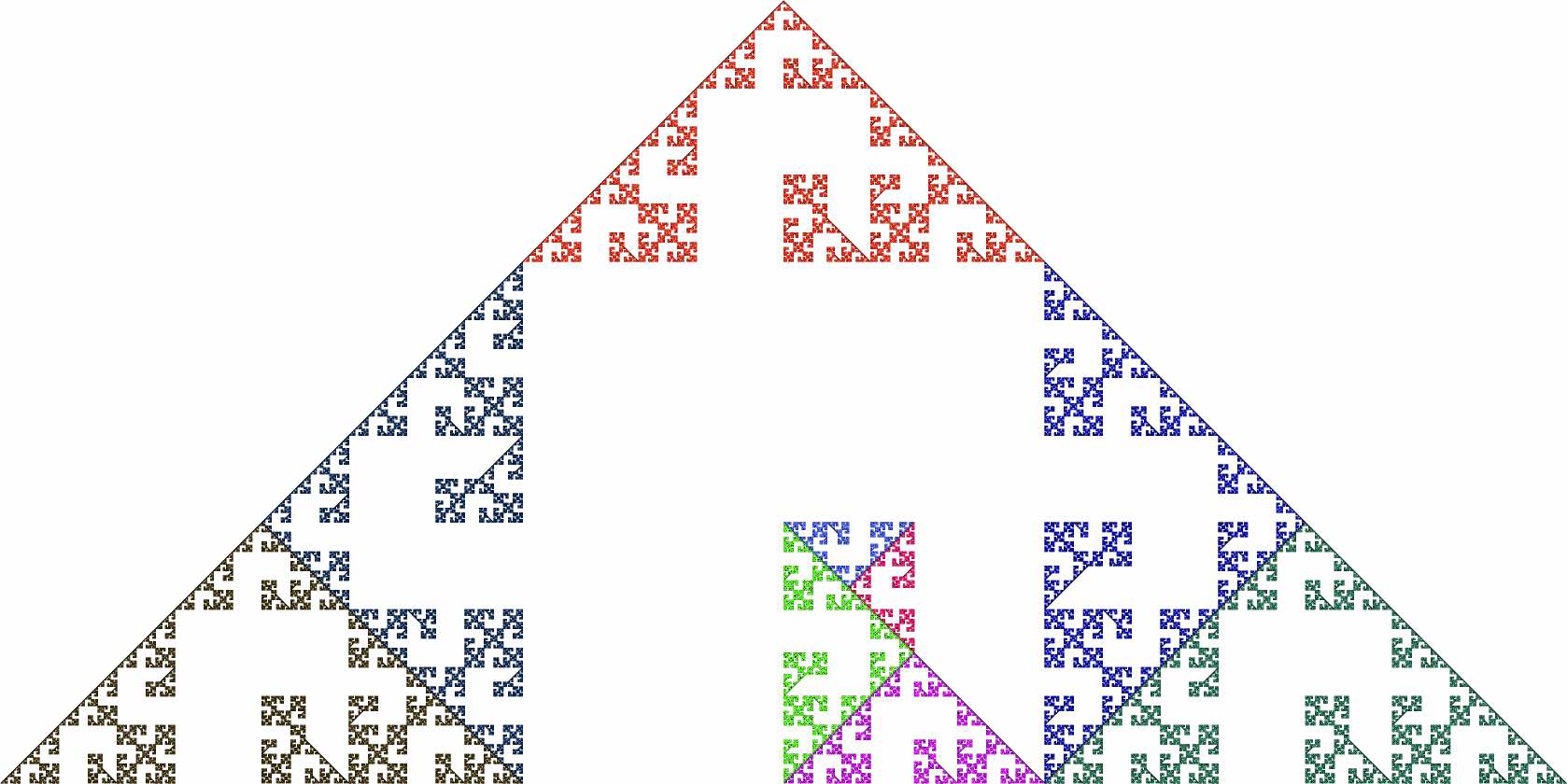}}; 
\draw (0.13,-0.25) node {$P_1 $} (0.18,-0.1) node  {$P_2 $} 
(0.37,0) node {$P_3 $}(0.58,-0.1) node  {$P_4 $} (0.64,-0.25) node {$P_5 $} (0.44,-0.25) node {$P_6 $} (0.35,-0.17) node {$P_7 $} (0.41,-0.07) node {$P_8 $} (0.47,-0.16) node {$P_9 $};
\draw (.94,-0.25) node {$K_1 $} (0.98,0) node  {$K_2 $} 
(1.18,0) node {$K_3 $}(1.39,0) node  {$K_4 $} (1.45,-.25) 
node {$K_5 $} (1.26,-0.25) node {$K_6 $} (1.155,-0.17) node {$K_7 $} (1.22,-0.07) node {$K_8 $} (1.275,-0.14) node {$K_9 $};
\end{tikzpicture}}}
 
 \begin{cor}
 Let $\eS$ be a generalized $P$-polygonal system, satisfying the condition(\ref{icnd}).
 For any subarc $\ga_{xy}\IN K$ and for any $n$, there is unique chain of pairwise different multiindices $\bi_1,\bi_2,...,\bi_l\in I^n$, which divides $\ga_{xy}$ to sequential arcs $\ga_{xx_1}\IN K_{i_1},\ldots,\ \ga_{x_{k-1}x_k}\IN K_{i_k},\ldots,\ \ga_{x_{l-1}y}\IN K_{i_l}$.\vse
 \end{cor}

\subsection{ $\da$-deformations of contractible polygonal systems.}

\begin{dfn}\label{deform} Ler $\da>0$.
A generalized $P'$-polygonal system $\eS'=\{S'_1,...,S'_m\}$ is called a $\da$-deformation of a $P$-polygonal system $\eS=\{S_1,...,S_m\}$, if there is a bijection $f:\bigcup\limits_{k=1}^m \eV_{P_k}\to \bigcup\limits_{k=1}^m \eV_{P'_k}$, such that\\
a) $f|_{\eV_P}$ extends to a homeomorphism $\tilde f: P\to  P'$; \\ b) $|f(x)-x|<\delta$  for any $x\in \bigcup\limits_{k=1}^m \eV_{P_k}$\\  c) $f(S_k(x))=S'_k(f(x))$ for any $k\in I$ and $x\in \eV_P$.\bigskip
\end{dfn}
\hspace{-1em}\resizebox{.9    \textwidth}{!}{{
\begin{tikzpicture}[x=10.0cm,y=10.0cm]
\node[anchor=south west,inner sep=0] at (0,0) {\includegraphics[width=.43\textwidth]{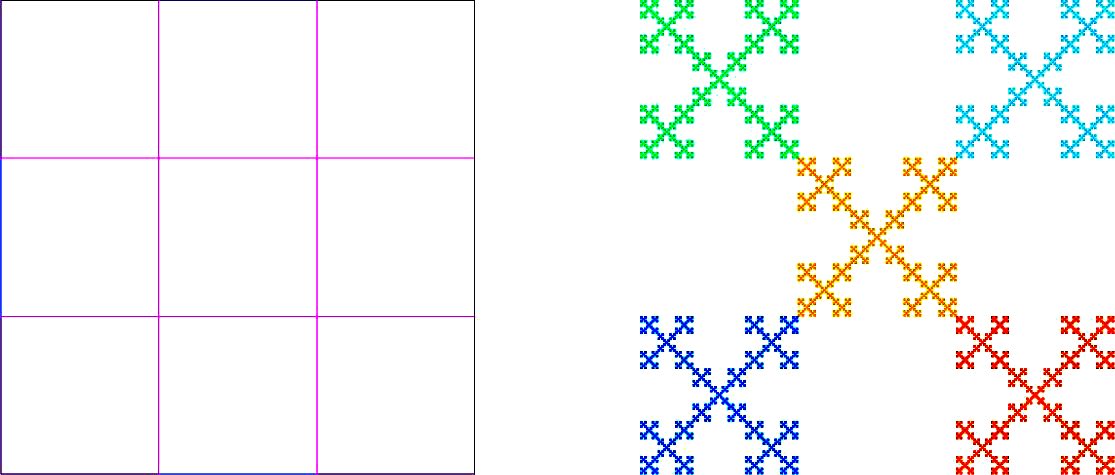}};
 \node[anchor=south west,inner sep=0] at (0.85,-.03){\includegraphics[width=.47\textwidth]{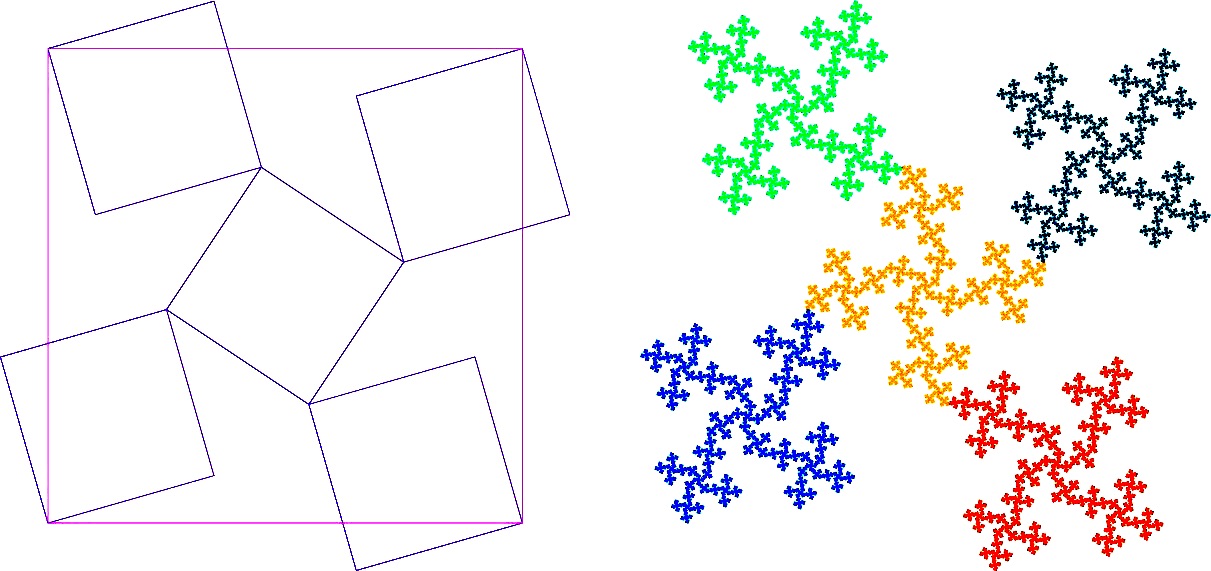}};
 \node[anchor=south west,inner sep=0] at (0.4,-.08){\small A polygonal system $\eS$ and its $\da$-deformation $\eS'$  };
    \node at(0.27,.25){$P_3$}; \node at (.15,.15) {$P_5$};\node at (.27,.05) {$P_2$};\node at (.05,.25) {$P_4$};\node at (0.05,.05) {$P_1$};
\node at (1.15,.24){$P'_3$}; \node at (1.05,.15) {$P'_5$};\node at (1.12,0.03) {$P'_2$};\node at (0.94,.27) {$P'_4$};\node at (.93,.06){$P'_1$};

\end{tikzpicture}}}\\

 Notice that by the Definition \ref{deform} if $z_1,z_2\in \eV_P$, $i,j\in I$  and $S_i(z_1)=S_j(z_2)$, then $S'_i( f(z_1))=S'_j(f(z_2))$. Moreover, we have the following
 \begin{lem}\label{bibj}
 If $A_1,A_2\in \eV_P$, $\bi,\bj\in\ia$ and $S_\bi(A_1)=S_\bj(A_2)$, then $S'_\bi( f(A_1))=S'_\bj( f(A_2))$.
 \end{lem}
 
\dok
Suppose $S_\bi(A)=B\in \eV_\wP$ for some $A\in \eV_P$ and let $\bi=i_1i_2...i_n$. Denote $S_{i_{k+1}...i_n}(A)$ by $A_k$.

Then we have a finite sequence of relations between $B\in \eV_\wP$ the vertices $A_k\in\eV_P$:\\
\beq\label{chn1} B=S_{i_1}(A_1); \quad A_1=S_{i_2}(A_2); \quad \ldots A_{n-1}=S_{i_n}(A)\eeq
Since, by c), $f(S_k(A_k))= S'_k(A'_k)$, \qquad $A'_{k-1}=f(A_{k-1})=f(S_k(A_k))=S'_k(A'_k)$, therefore the map $f$ transforms the relations \ref{chn1} to  
\beq\label{chn2} B'=S'_{i_1}(A'_1); \quad A'_1=S'_{i_2}(A'_2); \quad \ldots  A'_{n-1}=S'_{i_n}(A')\eeq
which implies $S'_\bi(A')=B'$\\
Therefore if $S_\bi(A_1)=S_\bj(A_2)\in \eV_{\wP}$, then $S'_\bi(f(A_1))=S'_\bj(f(A_2))$. \\

Now suppose $S_\bi(A_1)=S_\bj(A_2)$ and   $\bi=\bl\bi'$, $\bj=\bl\bj'$ and $S_\bi(A_1)=S_\bj(A_2)=S_\bl(B)$ for some $B\in \eV_\wP$. Then  $S_{\bi'}(A_1)=S_{\bj'}(A_2)=B$, therefore $S'_{\bi'}(f(A_1))=S'_{\bj'}(f(A_2))=f(B)$ and
$S'_{\bi}(f(A_1))=S'_{\bj}(f(A_2))=S'_\bl(f(B))$. \vse\\

\begin{thm}\label{attrmap}
Let $K$ and $K'$ be the attractors of a contractible $P$-polygonal system $\eS$ and of its $\da$-deformation $\eS'$ respectively and $\pi:I^\8\to K,\pi':I^\8\to K'$ be respective address maps.\\ (i)  There is unique continuous map $\hat f:K\to K'$ such that
$\hat f\circ\pi=\pi'$.\\ (ii) If $\eS'$ satisfies condition \ref{icnd}, then $\hat f$ is a homeomorphism.
\end{thm}
\begin{rmk}
Equivalent formulation of the statement (i) of the Theorem is:\\{\em There is unique continuous map $\hat f:K\to K'$ such that
 for any $z\in K$ and $\bi\in\ia$,}
 \beq\label{compat}\hat f(S_\bi(z))=S'_\bi(\hat f(z)).\eeq
\end{rmk}

\dok
The proof is similar to (cf.\cite[Lemma 1.]{ATK}). First, we define the function   $\hat f$ which is a surjection  of the dense subset  $G_\eS(\eV_P)\IN K$ to the dense subset $G_{\eS'}(\eV_{P'})\IN K'$.  Second, we show that it is H\"older continuous on $G_\eS(\eV_P)$, and therefore has unique continuous extension to a surjection of $K$ to $K'$, which we denote by the same symbol $\hat f$. Third, we show that the condition \ref{icnd} implies that $\hat f$ is injective and therefore is a homeomorphism.\\

1. Define a map $\hat f(z):G_\eS(\eV_P)\to G_{\eS'}(\eV_{P'})$ 
 by:
\beq\label{hatf}\hat f(z)=S'_\bi(f(S_\bi^{-1}(z))\mbox{  if   }z\in S_\bi(\eV_P)\eeq
As it follows from  Lemma \ref{bibj},
if $S_\bi(A_1)=S_\bj(A_2)=z$ then $S'_\bi(f(S_\bi^{-1}(z)))=S'_\bj(f(S_\bj^{-1}(z)))$, so the map $\hat f$ is well-defined.\\ Obviously, $\hat f(G_\eS(\eV_P))= G_{\eS'}(\eV_{P'})$  because if $A'\in \eV_{P'}$ and $z'=S'_\bi(A')$, then there is  a vertex $A=f^{-1}(A')\in \eV_P$, therefore $z'=\hat f(S_\bi(A))$.\\

Moreover, for any  $z\in G_\eS(\eV_P)$ and $\bi\in\ia$, $\hat f(S_\bi(z))=S'_\bi(\hat f(z))$ and\\
if $z_1,z_2\in G_\eS(\eV_P)$, $\bi,\bj\in\ia$ and $S_\bi(z_1)=S_\bj(z_2)$, then $S'_\bi(\hat f(z_1))=S'_\bj(\hat f(z_2))$.\\

2. Let $q_k = \Lip S_k$, $q_k' = \Lip S_k'$,  $\beta = \min \limits _{k\in I } {\dfrac {\log {q'_k}} {\log {q_k}}}$.\\

   Then, following the proof of \cite[Theorem 27, step 4.]{TSV0}, 
 in which for our estimates we use $K'$ instead of $P'$, 
   we see that for any $z_1,z_2\in G_\eS(\eV_P)$, $$|z_1'-z_2'| \le     \dfrac {2|K'|} {(\rho_0\cdot\sin {(\al_0/2)})^{\beta}}|z_1-z_2|^{\beta}.$$
   
 Therefore the map $\hat f$ can be extended to a H\"older continuous surjective map of $K$ to $K'$. Since for any $z\in K$ and any $k\in I$, $\hat f(S_k(z))=S'_k(f(z))$, $\hat f\circ\pi=\pi'$.\\
 
 3. Now, suppose the system $\eS'$ satisfies the condition (\ref{icnd}). Suppose for some $\bm\sa=i_1i_2...\in I^\8$ and $\bm\tau=j_1j_2...\in I^\8$, $\hat f\circ\pi(\bm\sa)=\hat f\circ\pi(\bm\tau)$. Then, if $i_1\neq j_1$, then, by condition \ref{icnd}, $P'_{i_1}\cap P'_{j_1}\neq\0$, therefore   $P_{i_1}\cap P_{j_1}=\{B\}$ for some $B\in \eV_\wP$ and $\pi(\bm\sa)=\pi(\bm\tau)=B$.
 
 Suppose now  $\bm\sa=\bl\bm\sa'$ and $\bm\tau=\bl\bm\tau'$ and $\hat f\circ\pi(\bm\sa)=\hat f\circ\pi(\bm\tau)$. Then, by formula
 \ref{compat}, $\hat f\circ\pi(\bm\sa')=\hat f\circ\pi(\bm\tau')$, so if first indices in $\bm\sa'$ and  $\bm\tau'$ are different, then $\pi(\bm\sa)=\pi(\bm\tau)=S_\bl(B)$ for some $B\in \eV_\wP$. 
 
 This implies injectivity of the map $\hat f$. So $\hat f$ is a homeomorphism of compact sets $K$ and $K'$. \vse\\

\subsection{Parameter matching theorem.}

The Definition \ref{cyclic} of cyclic vertices  can be applied to generalized polygonal systems. In this case, if $A$ is a cyclic vertex of a generalized P-polygonal system $\eS$, the map $S_\bi$ for which $S_\bi(A)=A$, need not be a homothety and we have to define the rotation parameter  for such map. Though the rotation angle $\al_i$ of the map $S_\bi$ is defined up to $2n\pi$, the number $n$ is unique defined by the set $\wP$.

\begin{dfn}
 Let $A$ be a cyclic vertex and $S_\bi(z)=r e^{i \al} (z-A)+A$, then the parameter $\lambda _A$ of the cyclic vertex $A$ is a number $\dfrac{\al}{\ln{r}}$, where the angle $\al$ is defined by the geometrical configuration of the system.
\end{dfn}

\begin{rmk} The following picture shows how the angle $\al$ depends on the geometric configuration of the system $\eS$.\end{rmk}
  \hspace{-1em}\resizebox{.9    \textwidth}{!}{{\small\begin{tikzpicture}[line cap=round,line join=round,>=stealth ,x=10.0cm,y=10.0cm]
\node[anchor=south west,inner sep=0] at (0,-.227) {
   \includegraphics[width=.24\textwidth]{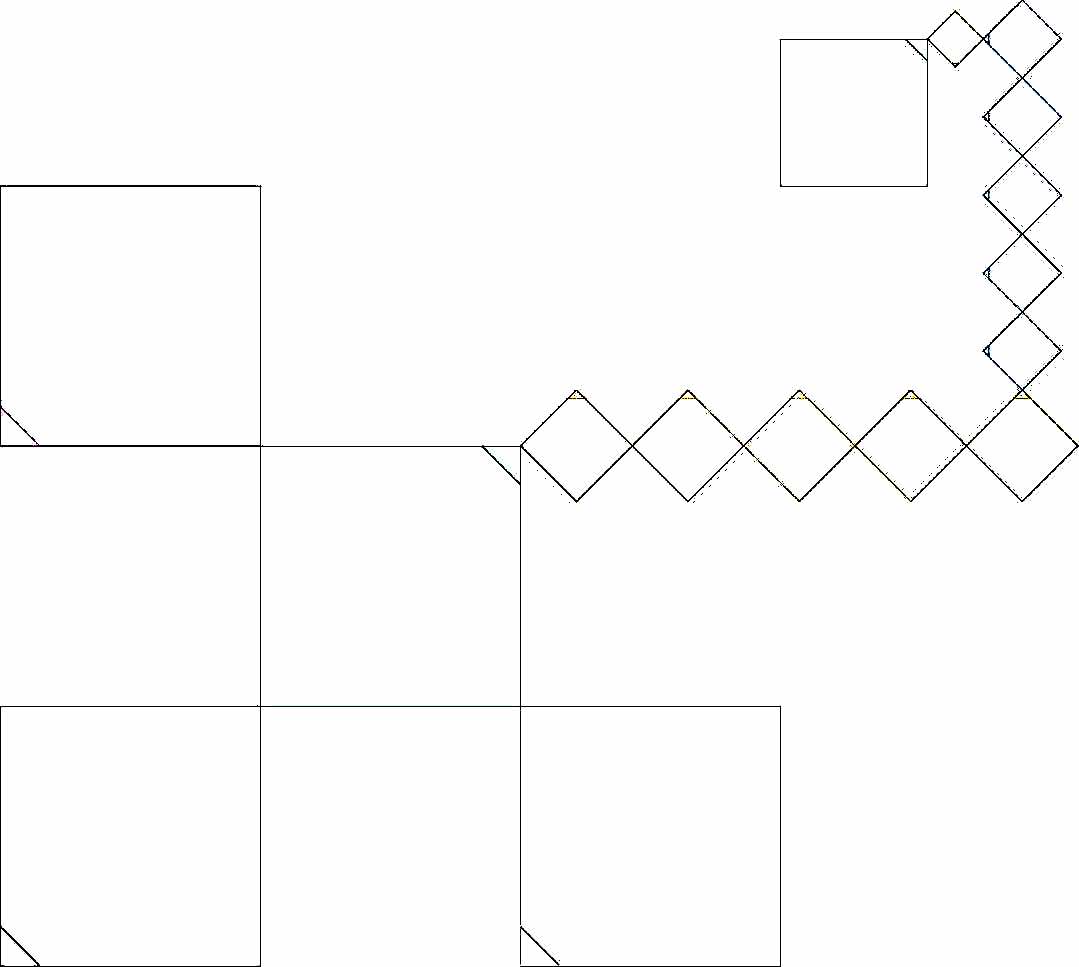}\quad  \includegraphics[width=.23\textwidth]{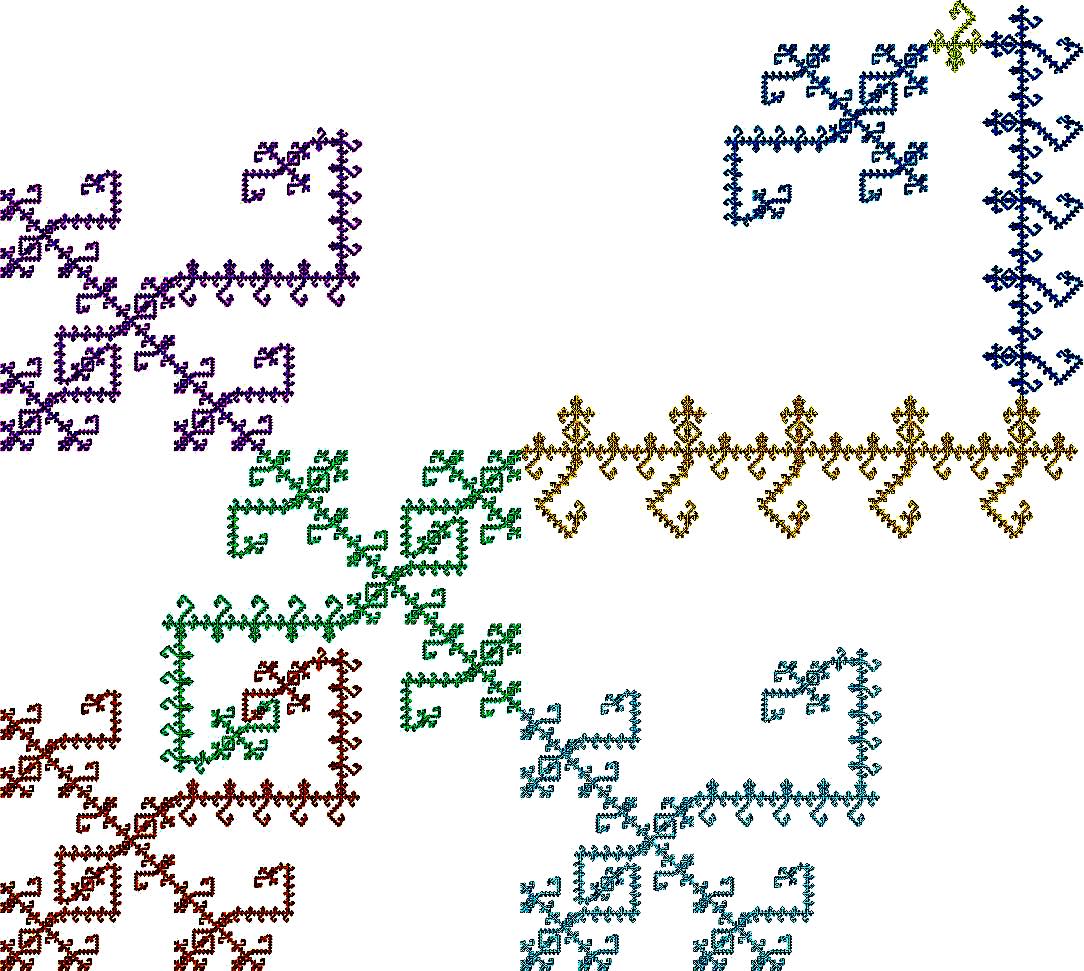}\quad
    \includegraphics[width=.22\textwidth]{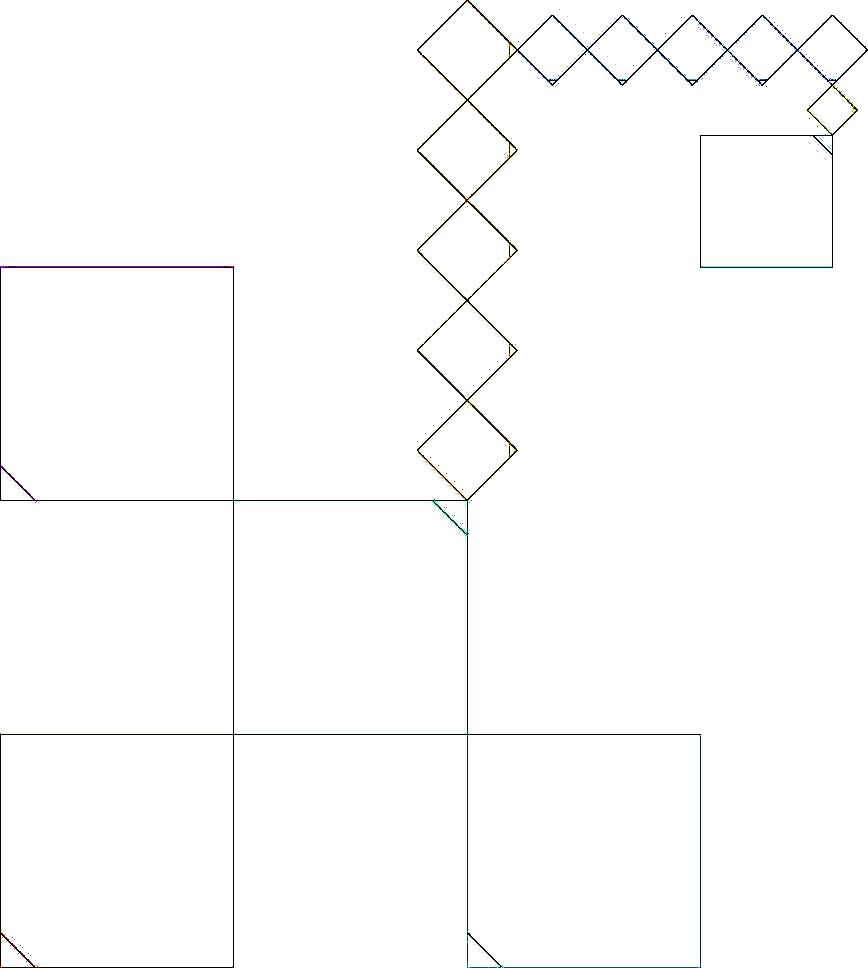}\quad  \includegraphics[width=.22\textwidth]{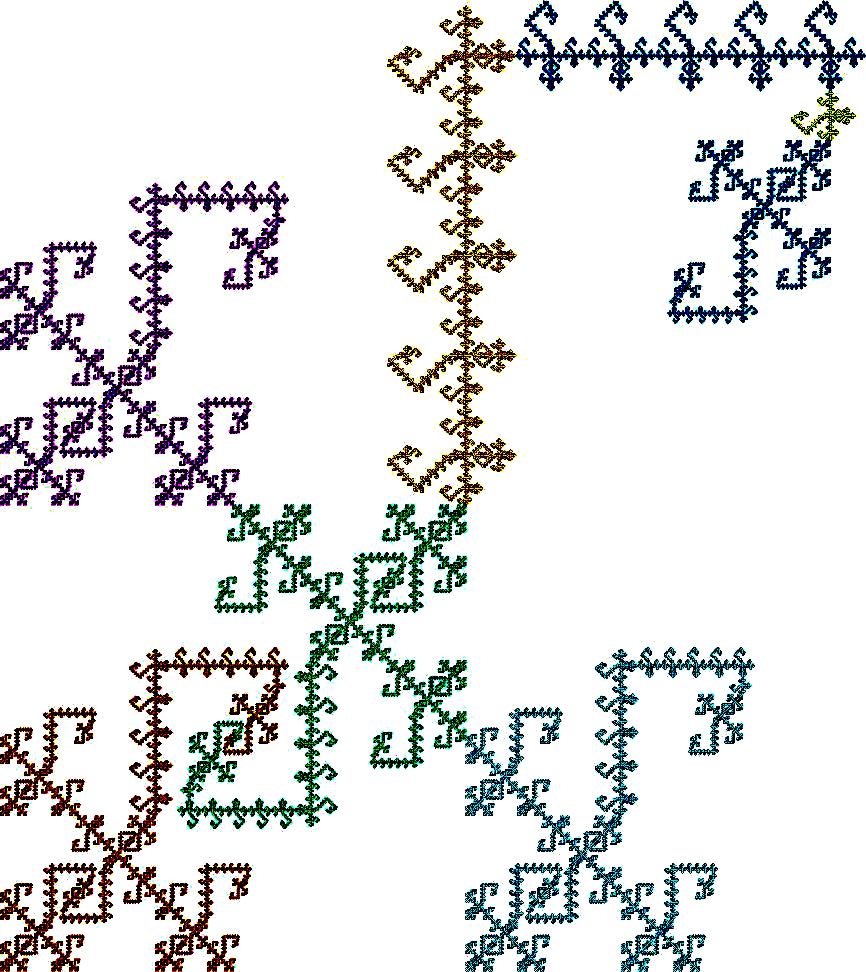}}; 
\draw (0.31,0.05) node {$\al=\pi$} (1.19,0.05) node {$\al=-\pi$};
\end{tikzpicture}}}

\begin{dfn} Generalized $P$-polygonal system $\eS$ of similarities satisfies the {\em parameter matching condition}, if for any $B \in\cup_{i=1}^m  \eV_{P_i}$ and for any cyclic vertices $A,A'$ such that for some $\bi,\bj\in I^*$,  $S_\bi(A)=S_\bj(A')=B$, the equality $\lambda _{A}=\lambda _{A'}$ holds.\end{dfn}

\begin{lem}
Let $\eS$ be a generalized $P$-polygonal system,  satisfying the condition(\ref{icnd}). For any vertices $A,B\in \eV_P$ there are $A', B'\in \eV_P$ and a map $S_i\in \eS$ such that $S_i(A')=A$ and $S_i(\ga_{A'B'})\IN\ga_{AB}$.
\end{lem}

\dok Consider the unique arc $\ga_{AB}$, connecting $A$ and $B$.

For the arc $\ga_{AB}$ we consider the chain $i_1,i_2,...,i_l$,  which partitions it to subarcs  $\ga_{Ax_1}\IN K_{i_1}$, ...,$\ga_{x_{k-1}x_k}\IN K_{i_k}$,...$\ga_{x_{l-1}B}\IN K_{i_l}$ (possibly to the only arc $\ga_{AB}$ if $\ga_{AB}\IN K_{i_1}$).  Put $A'=S_{i_1}^{-1}(A)$, $B'= S_{i_1}^{-1}(x_1)$, and $\ga(A'B')=S_{i_1}^{-1}(\ga_{Ax_1}).$\vse

\begin{prop}\label{fixparc}
  Let $\eS$ be a generalized $P$-polygonal system satisfying the condition (\ref{icnd}) and let $A$ be a cyclic vertex of the polygon $P$. Then there is such vertex $B\in V_{P}$ and a multiindex $\bi\in \ia$, that  $S_{\bi}(A)=A$  and the Jordan arc $\ga_{AB}\IN K$ satisfies the inclusion $S_{\bi}(\ga_{AB})\IN\ga_{AB}$.
\end{prop}

\dok 
Notice that if $\eS$ is a contractible $P$-polygonal system then for any cyclic vertex $A$ and for any $n$ there is {\em unique} multiindex $\bi\in I^n$, and unique vertex $B\in V_P$, such that $S_\bi(B)=A$. Therefore, if $S_\bi(A)=A$, the piece $S_\bi(K)$ separates the point $A$ from the other part of the attractor $K$ of the system  $\eS$, i.e. $A\notin \overline{K\mmm S_\bi(K)}$ and each Jordan arc $\ga_{AB}$ where $B\in V_P\mmm\{A\}$, contains a point $B'\in S_\bi(V_P\mmm \{A\})$. %

 In the case when $\eS$ is a generalized polygonal system, the situation is more complicated. Since the attractor $K$ is a  dendrite in the plane which has one-point intersection property, it follows from 
 \cite{TMV} that the system $\eS$ satisfies OSC and each vertex $A'\in V_P$ has finite ramification order. Let $U_1,...,U_s$ be the components of $K\mmm\{A\}$. Since $S_\bi$ fixes $A$, there is a permutation $\sa$ of the set $\{1,...,s\}$, such that for any $l\in\{1,...,s\}$, $S_\bi(U_l)\IN U_{\sa(l)}$. Therefore there is such $N$ that $\sa^N=\Id$ and $S_\bj=S_\bi^N$ sends each $U_l$ to a subset of $U_l$. Each of those components  $U_l$ which have non-empty intersection with $V_P\mmm\{A\}$ has also non-empty intersection with $S_\bj(V_P\mmm \{A\})$, therefore each arc  $\ga_{AB}, B\in V_P$ contains a point $B'\in S_\bj(V_P)$. \\

Let us enumerate the vertices of $P$  so that $A=A_1$ and other vertices are $A_2,...,A_{n_P}$. For each vertex $A_k, k\ge 2$ there is unique vertex $A_{k'}$ such that $\ga_{A_1A_k}\cap S_\bj(V_P)=S_\bj(A_{k'})$. The formula $\phi(k)=k'$ defines a map $\phi$ of $\{2,3...,n_P\}$ to itself. There is some $N'$ such that $\phi^{N'}$  has a fixed point $k$. Therefore $S_\bj^{N'}(\ga_{A_1A_k})\IN \ga_{A_1A_k}$.
\vse

\begin{dfn} 
The arc $\ga_{AB}$ is called an {\em invariant arc} of the cyclic vertex $A$.\\
\end{dfn}

From  Propositions \ref{1ordcyc} and \ref{fixparc} and V.V.Aseev's Lemma about  disjoint periodic arcs \cite[Lemma 3.1]{ATK} we come to the following Parameter Matching Theorem:

\begin{thm}\label{PMT}
Let the generalized $P'$-polygonal system $\eS'$ be a  $\da$-deformation of a contractible $P$-polygonal system  $\eS$ and the  attractor $K'$ of the system $\eS'$ be a dendrite. Then the system $\eS'$ satisfies parameter matching condition.
\end{thm}

\dok Let $\eS$ be a generalized polygonal system whose attractor $K$ is a dendrite. Let $C \in\cup_{i=1}^m  \eV_{P_i}$ and $A,A'\in \eV_P$ be such cyclic vertices that for some
$i,j\in I$, $S_i(A)=S_j(A')=C$. Denote the images $S_i(K)$ and $S_j(K)$ by $K_i,\ K_j$ respectively. Without loss of generality we can suppose that the point $C$ has coordinate $0$ in $\bbc$.
Since for some $\bi,\bj\in I^*$, $S_\bi(A)=A$ and $S_\bj(A')=A'$,
the maps $S_{b1}=S_iS_\bi S_i^{-1}$ and $S_{b2}=S_jS_\bj S_j^{-1}$ have $C$ as their fixed point and $S_{b1}(K_i)\IN K_i$ and $S_{b2}(K_j)\IN K_j$. Let $S_{b1}(z)=q_\bi e^{i\al_\bi}z$ and
$S_{b2}(z)=q_\bj e^{i\al_\bj}z$. So the parameters of the vertices $A$ and $A'$ will be $\la_1=\dfrac{\al_\bi}{\log q_\bi}$ and  $\la_2=\dfrac{\al_\bj}{\log q_\bj}$. Let $\ga_{AB}\IN K$ and $\ga_{A'B'}\IN K$
be invariant arcs for the vertices $A$ and $A'$. Let also $\ga_1=S_i(\ga_{AB})$ and $\ga_2=S_j(\ga_{A'B'})$. Then 
$S_{b1}(\ga_1)\IN \ga_1$ and $S_{b2}(\ga_2)\IN \ga_2$. By V.V.Aseev's Lemma  on disjoint periodic arcs \cite[Lemma 3.1]{ATK} it follows that if $\ga_1\cap\ga_2=\{C\}$, then $\la_1=\la_2$.\vse\\

\subsection{Main theorem.}

{\bf Some assumptions.}   From now on we will use the following convention: $\eS=\{S_1,...,S_m\}$ will denote  a contractible $P$-polygonal system and $\eS'=\{S'_1,...,S'_m\}$ will be a $P'$-polygonal system which is a $\da$-deformation of $\eS$  defined by a map $f$.\\ For any $k\in I$, $S_k(z)=q_k e^{i\al_k}(z-z_k)+z_k$ and $S'_k(z)=q'_k e^{i\al'_k}(z-z'_k)+z'_k$, where $z_k=\fix(S_k)$. We also suppose by default that $\diam P=1$. We   suppose that \beq\label{assum}\da<q_{min}/8\mbox{\quad and \quad}\da<(1-q_{max})/8\eeq

\begin{lem}\label{lasin}
Let  $\eS'=\{S'_1,...,S'_m\}$ be a  $\da$-deformation of a contractible $P$-polygonal system $\eS$. For sufficiently small $\da$, and for any $k\in I$,
\beq\label{asin} \dfrac{q_k-2\da}{1+2\da}\le q'_k\le \dfrac{q_k+2\da}{1-2\da}  \mbox{\quad \rm and  \quad}|\al_k'-\al_k|\le \arcsin 2\da+\arcsin \dfrac{2\da}{q_k}.\eeq
\end{lem} 
\dok Let $A,B$ be such vertices of $P$ that $|B-A|=1$. Let us write $S_k(A)=A_k$ and $f(A)=A'$ and use the similar notation for all vertices so by definition, $S'_k(A')=A'_k=f(A_k)$. Notice that $\dfrac{B_k-A_k}{B-A}=q_k e^{i\al_k}$ and $\dfrac{B'_k-A'_k}{B'-A'}=q'_k e^{i\al'_k}$.

Since the map $f$ moves $A,B,A_k,B_k$ to a distance $\le\da$, so   $|(B-A)-(B'-A')|\le 2\da$ and $|(B_k-A_k)-(B'_k-A'_k)|\le 2\da$. Therefore 
$|(B_k-A_k)|-2\da\le|(B'_k-A'_k)|\le |(B_k-A_k)|+2\da$ and
\beq\label{dal}\al'_i-\al_i=\arg\dfrac{B'_k-A_k'}{B'-A'}\dfrac{B-A}{B_k-A_k}=\arg\dfrac{B'_k-A_k'}{B_k-A_k}-\arg\dfrac{B'-A'}{B-A}\eeq
This implies the inequalities (\ref{asin}).  \vse\medskip

Under the Assumptions (\ref{assum}), $3q_{min}/5<\dfrac{q_{min}-2\da}{1+2\da}<q'_k<\dfrac{q_{max}+2\da}{1-2\da}<\dfrac{1+3q_{max}}{3+q_{max}}$;\\ 
taking into account that $q_k<1$ and $1-2\da>3/4$, and that
if $0<x<.5$, then $\arcsin x<1.05x$, we have
 \beq\label{dadq}\Da q_k=|q'_k-q_k|<\dfrac{2\da(1+q_k)}{1-2\da}<6\da  \mbox{\quad and \quad} \Da \al_k=|\al'_k-\al_k|< 
C_\al\da\eeq where $C_\al=2.1(1+1/q_{min})$.\smallskip\\

Let $V_\da(P)$ denote $\da$-neighborhood of the polygon $P$.

\begin{lem}
Let  $\eS'=\{S'_1,...,S'_m\}$ be a  $\da$-deformation of a contractible $P$-polygonal system $\eS$. The set $U=V_{\da_1}(P)$, where $\da_1=\dfrac{8\da}{1+3q_{max}}$, satisfies the condition
\beq \label{vda1}\mbox{\rm  for any  } k\in I,\ \ \ \    S_k(U)\IN U\mbox{\quad\rm and \quad}S'_k(U)\IN U\eeq\end{lem} 

\dok
By Definition \ref{deform}, $V_\da(P_k)\NI P'_k$, $V_\da(P'_k)\NI P_k$ and since vertices of $P$ are also moved at distance less than $\da$,
$V_\da(P)\NI P'$ and $V_\da(P')\NI P$.

 So we can write $S_k'(P')\IN V_\da(P_k)\IN V_\da(P)$ from which it follows that
$S_k'(P)\IN V_{2\da}(P_k)\IN V_{2\da}(P)$.\\ For any positive $\rho$ we have the inclusion 
$S_k'(V_\rho(P))\IN V_{2\da+q'_k\rho}(P)$. In the case when $\rho= 2\da+q'_k\rho$ this implies $S_k'(V_{\rho}(P))\IN V_{\rho}(P)$ where $\rho=\dfrac{2\da}{1-q_k'}$. Since $q_k'\le 
q_k+2\da$, $q'_{max}\le q_{max} +2\da<\dfrac{3q_{max}+1}{4}$, we come to inclusions (\ref{vda1}).\vse

\begin{lem}
For any $z\in V_{\da_1}(P)$, $|S'_k(z)-S_k(z)|<C_\Da\da$, where $C_\Da=14+2C_\al$.
\end{lem}
\dok Take $z\in V_{\da_1}(P)$ and consider the difference $S'_k(z)-S_k(z)$. It can be represented in the form $S'_k(A)-S_k(A)+(q'_ke^{i\al'_k}-q_ke^{i\al_k})(z-A)$. So \beq\label{delta}|S'_k(z)-S_k(z)|<|S'_k(A)-S_k(A)|+(|q'_k-q_k|+q_k|e^{i\al'_k}-e^{i\al_k}|)|z-A|.\eeq 
Since $|z-A|<1+\da_1<2$ and $|S'_k(A)-S_k(A)|<2\da$, the right hand side of (\ref{delta}) is no greater than $2\da+2(6\da+C_\al\da)$.\vse

\begin{prop}\label{deltaK}
Let $\pi:I^\8\to K$ and $\pi':I^\8\to K'$ be the address maps for the systems $\eS$ and $\eS'$ respectively.\\
1.Under the assumptions (\ref{assum}), for any $\sa\in I^\8$, 
\beq\label{dKeq} |\pi'(\sa)-\pi(\sa)|< C_K\da \mbox{   where   }C_K=\dfrac{2C_\Da}{1-q_{max}}\eeq 
2. For any n, if  the system $\eS^{'(n)}$ is a generalized polygonal system, then it is  $C_K\da$- deformation of the system $\eS^{(n)}$.
\vse
\end{prop}

\begin{rmk}\label{Sshift}
Let  $\eS'=\{S'_1,...,S'_m\}$ be a  $\da$-deformation of a contractible $P$-polygonal system $\eS$. Let $A\in S_j(\eV_P)$ for some $j\in I$. Let $g(z)=z-A+A'$ and $\hat S''_k=g\circ S'_k\circ g^{-1}$.
Then $\eS''=\{S''_1,...,S''_m\}$ is a $2\da$-deformation of the system $\eS$, for which $A''=A$, $K''=g(K')$, $P''_\bj=g(P_\bj)$. Since $g$ is a translation, the estimates (\ref{asin}) and (\ref{dadq}) for $\eS''$ remain the same with the same $\da$, while $|\pi''(\sa)-\pi(\sa)|< (C_K+1)\da$.Thus we will denote $\da_2=(C_K+1)\da$. \end{rmk}

Taking into account the Propositions \ref{1ordcyc} and \ref{deltaK}, it is sufficient to prove the Theorem for the case when all cyclic vertices of the system $\eS$ have order 1.

\begin{prop}
Let $P'$-polygonal system $\eS'$ be a $\da$-deformation of a contractible $P$-polygonal system $\eS$. Let $A\in \eV_P$ be a cyclic vertex (of order 1) and $S_k(z)=q_ie^{i\al_k}(z-A)+A$.  Then the rotation angle $\al_k$ of the map $S'_k$ does not exceed $\arcsin 2\da+\arcsin \dfrac{2\da}{q_k}$ and the parameter $\la_k$ of the map $S'_k$  satisfies the inequality \begin{equation}\label{parmeq}|\la_k|\le\dfrac{\arcsin 2\da+\arcsin\dfrac{2\da}{q_k}}{|\log(q_k+2\da)-\log(1-2\da)|}\end{equation}

\end{prop}
\dok The formula (\ref{parmeq}) follows directly from  Lemma \ref{lasin}.\vse\\

Under the assumptions (\ref{assum}), \beq\label{prmeq2}|\la_k|<C_\la\da\mbox{, where }C_\la=\dfrac{2.1(1+1/q_{max})}{\log(3+q_{max})-\log(3q_{max}+1)}.\eeq

\begin{lem}
Let $\eS$ be a contractible $P$-polygonal system whose cyclic vertices have order 1 and $\eS'$ be its $\da$-deformation. Then if \beq \label{mineq}2.1\dfrac{\da_2}{\rho_1}+\la\log\dfrac{\rho_2+\da_2}{\rho_1-\da_2}<\al_0 \mbox{ and }  2\da_2<\rho_0,\eeq then the system $\eS'$ satisfies the Condition (\ref{icnd})
\end{lem}

{\bf Proof.}  Take a vertex $B\in V_\wP$. We may suppose for convenience that $B=0$ and, following Remark \ref{Sshift}, we can suppose that the mapping $f$ fixes the vertex $B=0$, so $B'=B=0$.   Let $W_l=S_{\bj_l}(K\mmm K_{i_l})$. The maps $\bar S_l=S_{\bj_l}S_{i_l}S_{\bj_l}^{-1}$ are homotheties with a fixed point $B$ such that 
\beq\label{kbj} K_{\bj_l}\mmm\{B\}=\bigsqcup\limits_{n=0}^\8\bar S_l^n(W_l)\eeq 
Similarly, let $W'_l=\hat f(W_l)$ and $\bar S'_l=S'_{\bj_l}S'_{i_l}S_{\bj_l}^{'-1}$. Then
 \beq \label{kpbj} K'_{\bj_l}\mmm\{B\}=\bigsqcup\limits_{n=0}^\8\bar S_l^{'n}(W'_l)\eeq
Notice that for any $l$, $\bar S_l(z)=q_{i_l}z$ and $\bar S'_l(z)=q'_{i_l}e^{i\al_{i_l}}z$, and due to parameter matching condition, there is such $\la$, that for any $l$, $\al_{i_l}=\la\log q'_{i_l}$. 
  
Consider the map $z = exp(w)$ of the plane $(w = \ro+i\fy)$ as universal cover of the punctured plane $\bbc\mmm\{0\}$. 

Consider   polygons $P_{\bj_l}$ and choose their liftings in the plane $(w = \ro+i\fy)$.
We may suppose these liftings lie in respective horizontal strips $\te^-_l\le\fy\le\te^+_l$, where $0<\te^-_l<\te^+_l<2\pi$ and $\te^+_l+\al_0<\te^-_{l+1}$ for any $l<k$ and $\te^+_k+\al_0<\te^-_1 +2\pi$.
We also consider liftings of $K_{\bj_l}$, $W_l$, $K'_{\bj_l}$ and
$W'_l$. We denote these liftings by $\eK_{\bj_l}$, $\eW_l$, $\eK'_{\bj_l}$ and
$\eW'_l$.
It follows from the equations \ref{kbj} and \ref{kpbj}, that
\beq\label{ekbj} \eK_{\bj_l}=\bigsqcup\limits_{n=0}^\8\bar T_l^n(\eW_l) \mbox{\quad and \quad} \eK'_{\bj_l}=\bigsqcup\limits_{n=0}^\8\bar T_l^{'n}(\eW'_l) \eeq  
where $T_l(w)=w+\log q_l$   and   $T'_l(w)=w+(1+i\la)\log q'_l$ are parallel translations for which $T_l(\eK_l)\IN \eK_l$ and $T'_l(\eK'_l)\IN \eK'_l$.

\resizebox{.95    \textwidth}{!}{{
\begin{tikzpicture}[line cap=round,line join=round,>=stealth ,scale=1]\scriptsize
\node at (-8.5,-2.5) {
   \includegraphics[width=.28\textwidth]{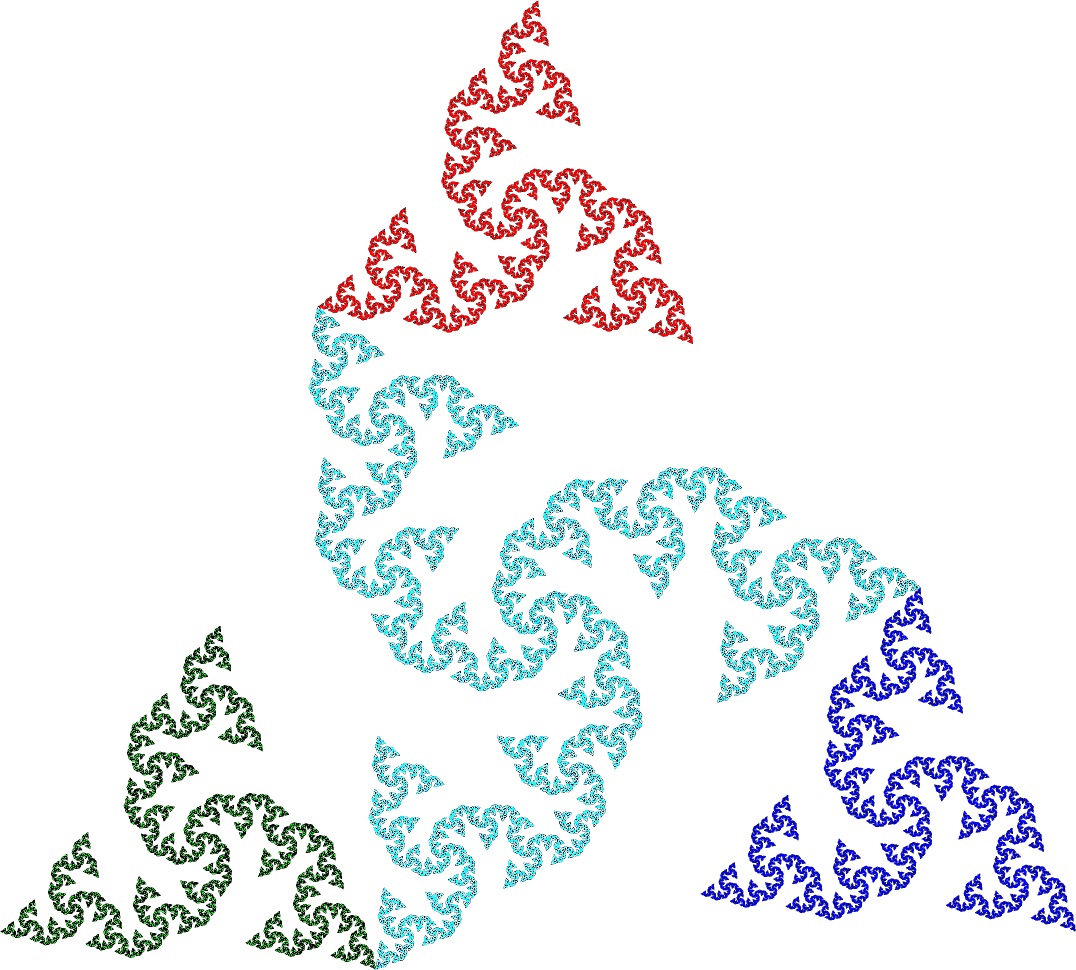}}; \node at (-2 ,-2.68) {\includegraphics[width=.35\textwidth]{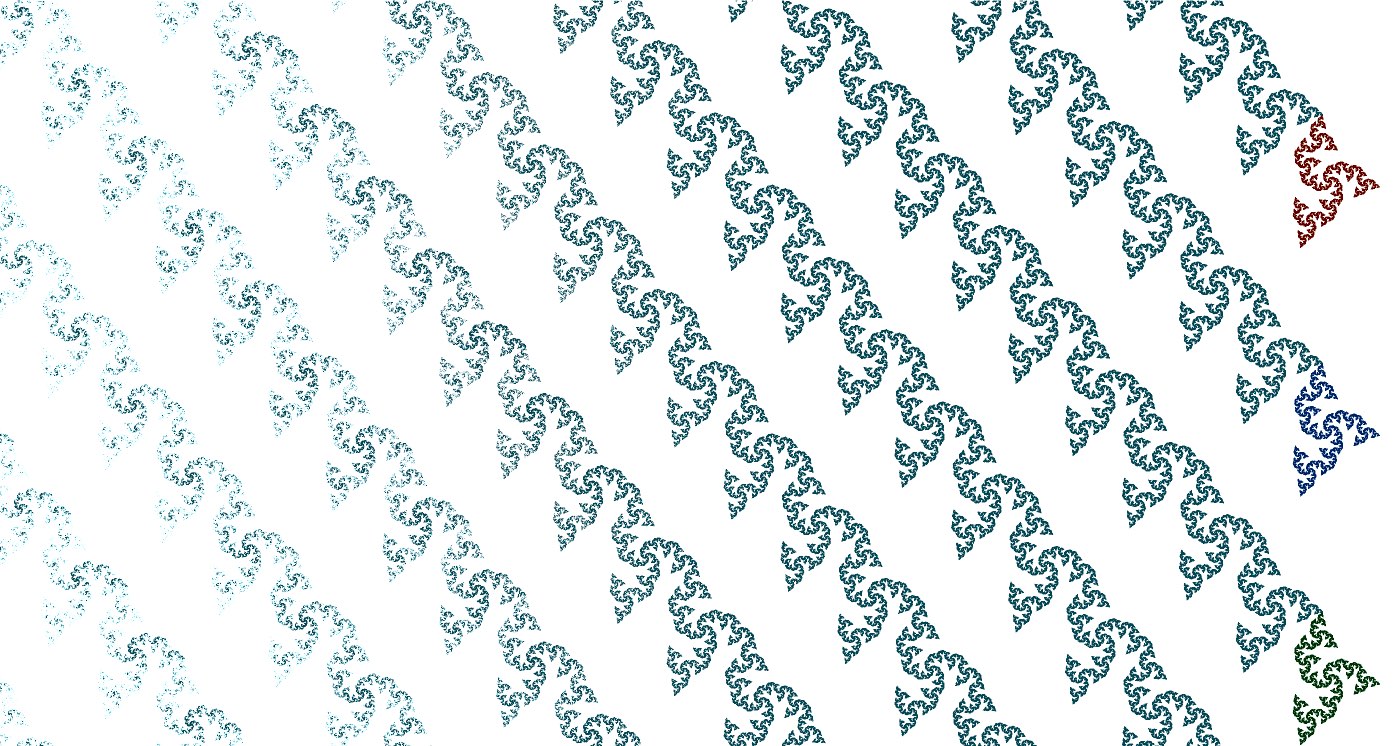}};\node at (4.8,-2.68) {\includegraphics[width=.35\textwidth]{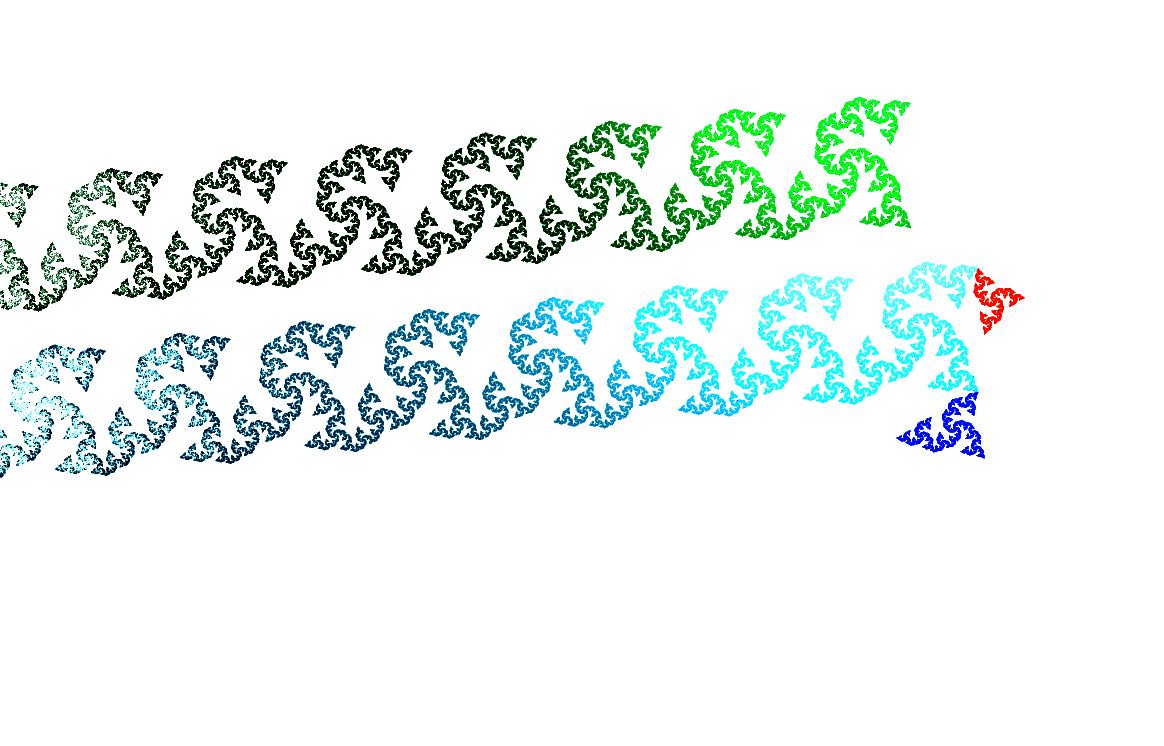}};
	    \node at(-9.2,-4.77){$B$};\node at(-8.65,-3.1){$O$};
	    \node at (-1.7 ,-5.3) {\small The images of the set $K'$ \qquad \qquad under   the map $w=\log(z-O)$  \qquad\qquad and   the map  $w=\log(z-B)$.};
\end{tikzpicture}}}

The sets $\eK_l$ lie in the half-strips $ \ro\le \log\rho_2, \te^-_l\le\fy\le\te^+_l$, while the sets $\eW_l$ are contained in rectangles $R_l=\{\log\rho_1 \le\ro\le \log\rho_2, \te^-_l\le\fy\le\te^+_l\}$.

Then the sets $\eW'_l$ lie in a rectangle $$R'_l=\left\{\log(\rho_1-\da_2)\le\ro\le \log(\rho_2+\da_2), \te^-_l-1.05\dfrac{\da_2}{\rho_1}\le\fy\le\te^+_l+1.05\dfrac{\da_2}{\rho_1}\right\}$$

Each union  $\bigcup\limits_{n=0}^\8T^{'n}_l(R'_l)$ lies in a half strip 
$$ \begin{cases}\ro\le \log(\rho_2+\da_2)\\   \te^-_l-1.05\dfrac{\da_2}{\rho_1}-\la \log(\rho_2+\da_2)\le\fy-\la\ro\le \te^+_l+1.05\dfrac{\da_2}{\rho_1}-\la \log(\rho_1-\da_2)\end{cases}$$

Therefore the set $\eK'_{\bj_l}$ also lies in this half-strip.
So, if 
\beq\te^+_{l-1}+1.05\dfrac{\da_2}{\rho_1}-\la \log(\rho_1-\da_2)<\te^-_l-1.05\dfrac{\da_2}{\rho_1}-\la \log(\rho_2+\da_2)\eeq then $\eK'_{\bj_{l-1}}\cap\eK'_{\bj_l}=\0$.

We can guarantee that such inequality holds for any $l$ if 
$2.1\dfrac{\da_2}{\rho_1}+\la\log\dfrac{\rho_2+\da_2}{\rho_1-\da_2}<\al_0$.

If, moreover, $2\da_2<\rho_0$, then for any $i_1,i_2\in I$  such that $P_{i_1}\cap P_{i_2}=\0$, $P'_{i_1}\cap P'_{i_2}=\0$  and $K'_{i_1}\cap K'_{i_2}=\0$ which implies the condition (\ref{icnd}).
 \vse

 \begin{thm}\label{mainthm}
Let $\eS$ be a contractible $P$-polygonal system. There is such $\da>0$ that for any $\da$-deformation $\eS'$ of the system $\eS$, satisfying parameter matching condition, the attractor $K(\eS')$ is a dendrite, homeomorphic to $K(\eS)$.
\end{thm}

\dok  Let all the cyclic vertices of the $P$-polygonal system $\eS$ have order 1. If we suppose that $\da_2<\rho_1/4$,and $\da_2<(1-\rho_2)/4$ and
combine the inequalities \ref{dadq},\ref{dKeq},\ref{prmeq2},\ref{mineq}, we see that if
the following inequalities hold:\\
1.$\da<\dfrac{q_{min}}{8}$;\quad  2. $\da<\dfrac{1-q_{max}}{8}$;\quad  3. $\da<\dfrac{\rho_0}{2(C_K+1)}$;\quad  4. $\da<\dfrac{\rho_1}{4(C_K+1)}$;\\

 5. $\da<\dfrac{1-\rho_2}{4(C_K+1)}$; \qquad and \qquad  6. $\da<\dfrac{\al_0}{\dfrac{2.1(C_K+1)}{\rho_1}+C_\la\log\dfrac{1+3\rho_2}{3\rho_1}}$,\\ 
 
 then the attractor $K'$ of $\da$-deformation $\eS'$ of the system $\eS$ satisfies  the condition (\ref{icnd}). Therefore $K'$ is a dendrite. By Theorem \ref{attrmap}, the map $\hat f:K\to K'$ is a bijection and therefore it is a homeomorphism.
 
 Suppose now that $\eS$ has cyclic vertices of order greater than 1 and let $M=12+4.2\left(1+\dfrac{1}{q_{min}}\right)$. There is such $n$, that the system $\eS^{(n)}$ has cyclic vertices of order 1. Suppose any $\da$-deformation of the system $\eS^{(n)}$ generates a dendrite. Then for any $\da/M$-deformation deformation $\eS'$ of the system $\eS$, the system $\eS^{'(n)}$ is a $\da$-deformation of the system $\eS^{(n)}$.\vse

\end{document}